\documentclass[reqno,round]{alea2}
\usepackage{natbib}
\usepackage{fancyhdr}
\usepackage{graphicx}

\usepackage{amssymb,amsfonts,chemarr}
\usepackage{amsmath}

\makeatletter \@addtoreset{equation}{section} \makeatother

\renewcommand\theequation{\thesection.\arabic{equation}}
\renewcommand\thefigure{\thesection.\@arabic\c@figure}
\renewcommand\thetable{\thesection.\@arabic\c@table}

\theoremstyle{plain}
\newtheorem{theorem}{Theorem}[section]                                          
\newtheorem{prop}[theorem]{Proposition}                          
\newtheorem{lem}[theorem]{Lemma}

\theoremstyle{definition}

\theoremstyle{remark}
\newtheorem{remark}[theorem]{Remark}

\def\prf{\noindent{\textit{Proof: }}}

\def\C{\mathbb{C}}
\def\H{\mathbb{H}}
\def\E{\mathbb{E}}
\def \b{{\tt b}}
\def \e{{\rm e}}
\def \d{{\tt d}}
\def \be{\begin{eqnarray*}}
\def \ee{\end{eqnarray*}}
\def \ben{\begin{eqnarray}}
\def \een{\end{eqnarray}}
\def\wh{\widehat}
\def\wt{\widetilde}
\def\noi{\noindent\underbar}
\def \sous#1#2{\mathrel{\mathop{\kern 0pt#1}\limits_{#2}}}
\def \sur#1#2{\mathrel{\mathop{\kern 0pt#1}\limits^{#2}}}
\def \el{\sur{=}{(d)}}

\begin{document}


\author{Alain Rouault}
\address{LMV  B\^atiment Fermat\\
 Universit\'e
Versailles-Saint-Quentin \\
F-78035 Versailles\\ France}
\email{Alain.Rouault@math.uvsq.fr}

\title[Asymptotic behavior of
random determinants]{Asymptotic behavior of random determinants in
the  Laguerre, Gram  and Jacobi ensembles}

\begin{abstract}
We consider
 properties of determinants of some random symmetric
matrices issued from multivariate statistics: Wishart/Laguerre
ensemble (sample covariance matrices), Uniform Gram ensemble (sample
correlation matrices) and Jacobi ensemble (MANOVA). If $n$ is the
size of the sample, $r\leq n$ the number of variates and $X_{n,r}$
such a matrix, a generalization of the Bartlett-type theorems gives
a decomposition of $\det X_{n,r}$ into a product of $r$ independent
gamma or beta random variables. For $n$ fixed, we study the
evolution as $r$ grows, and then take the limit of large $r$ and $n$
with $r/n = t \leq 1$. We derive limit theorems for the sequence of
{\it processes with independent increments} $\{n^{-1} \log \det X_{n
, \lfloor nt\rfloor}, t \in [0, T]\}_n$ for
 $T \leq 1$ : convergence in probability, invariance
principle, large deviations. Since the logarithm of the determinant
is a linear statistic of the empirical spectral distribution, we
connect the results for marginals (fixed $t$) with those obtained by
the spectral method. Actually, all the results hold true for log
gases or $\beta$ models, if we define the determinant as the product
of charges. The classical matrix models (real, complex, and
quaternionic) correspond to the particular values $\beta = 1,2,4$ of
the Dyson parameter.  
\end{abstract}
\maketitle

\section{Introduction}

Random determinants of symmetric matrices are of constant use in
random geometry to compute volumes of parallelotopes (see 
\citet{Nielsen1},  \citet{Mathai}) and in multivariate
statistics to build tests. Twenty years after the book of 
\citet{Girkobook}, recent developments in Random Matrix Theory add
a new interest to the study of their asymptotical behavior and invite to a new
insight.

Let $B= [b_1, \dots , b_r]$ be the $n\times r$ matrix  with $r$ column vectors $b_1, \dots , b_r$ of $\mathbb {R}^n$.
If $B'$ denotes its transpose, the determinant of the $r\times r$ Gram matrix $B'B$ 
satisfies the well known Hadamard inequality :
\ben
\label{Had} \det B'B \leq
\Vert b_1\Vert^2 \cdots\Vert b_r\Vert^2 \een with equality if and
only if  $b_1, \dots , b_r$ are orthogonal \citep{Hadamard}. It means that the
volume (or $r$-content) of the parallelotope built from  $b_1,
\dots, b_r$ is maximal when the vectors are orthogonal.
The quantity
\be h(B) =\frac{\det B'B}{\Vert b_1\Vert^2 \cdots\Vert b_r\Vert^2}
\,\ee is usually called the Hadamard ratio. If
we replace
sequentially $b_i$ by its projection $\widehat b_i$ on the orthogonal of the
subspace spanned by $b_1, \dots, b_{i-1}$ (Gram-Schmidt
orthogonalization), we have
\[\det B'B = \prod_{i=1}^r \Vert \widehat b_i\Vert^2\,.\]
Motivated by basis reduction
problems,  \citet{Schnorrortho} defined the orthogonality defect as the quantity
$1/\sqrt{h(B)}$  (see also  \citet{akhavi1} and references therein). \citet*{abbott}
and \citet{dixon} are concerned with the tightness of the bound
$h(B) \leq 1$ when $B$ is random and $n=r$. For these authors, the random vectors $b_i$ are sampled independently and uniformly on the
 unit sphere \[\mathbb{S}_{\mathbb R}^{n}  =\{(x_1,
\dots, x_n) \in \mathbb {R}^n :  x_1^2 + \dots + x_n^2 = 1 \}\,.\] It is known that then the variables $\Vert \widehat b_i\Vert^2$ are independent and
Beta distributed with varying parameters. When the entries of the matrix $B$ are independent and ${\mathcal N}(0, 1)$,  the variables $\Vert \widehat b_i\Vert^2$ are independent and
Gamma distributed with varying parameters \citep{Bart}.

Writing $B_{n,r}$
instead of $B$ to stress on dimensions and $X_{n,r} = B_{n,r}'B_{n,r}$, we are interested in this paper in the asymptotic behavior of
$\det X_{n,r}$ when $n$ and $r$ both tends to infinity, in the regime $r/n\rightarrow c\in [0,1]$.
Since the construction of the $\widehat b_i$ is recursive, it is possible (for fixed $n$) to consider the whole sequence of variables
$\{\det X_{n,r}, r=1, \dots,n\}$ at the same time.

It corresponds to
the decomposition of the determinant of a symmetric positive matrix
$A$ as \be \det A = \prod_{j=1}^r \frac{\det A^{[j]}}{\det A^{[j-1]}}\,, \ee
where $A^{[j]}$  is the $j\times j$  upper-left corner of $A$ with the convention $\det A^{[0]} = 1$.  When using this approach, we will refer to it as
the decomposition method.

The decomposition method is also valid when entries of the matrix
are complex, considering the Hermitian conjugate $B^\star$ and then
$B^\star B$, and also when the entries are real quaternions,
considering the dual $B^\dag$ and then $B^\dag B$.

In these three
cases, a Bartlett-type theorem gives the determinant as a product of
independent variables, with Gamma  or Beta distributions.
 Passing to
 logarithms, it is then possible to consider a triangular array of
variables and a   process with independent increments $\{n^{-1} \log \det X_{n
, \lfloor nt\rfloor}, t \in [0, T]\}_n$ for
 $T \leq 1$ indexed by the "time" $t= r/n$. Thanks to the additive structure of the log det, we obtained limit theorems : convergence in probability, invariance
principle and  large deviations.

 The same is
true for random matrices following the  Jacobi (or MANOVA)
distribution. Actually, the whole construction is possible in the
so-called $\beta$-models, which are an extension of the above ones,
which correspond to the three-fold way $\beta=1,2,4$ of Dyson. For other values of $\beta$ they are not defined as matrix models but log gases models, in which the eigenvalues are replaced by charges and determinants by  products of charges. It has be shown recently that they correspond also to models of
tri-diagonal random matrices  (see 
\citet*{EdelmanDumitriu}, \citet*{Killip1}, \citet*{edelman:bjm}).

Of course, for $r$ fixed, there is also another underlying structure of product:
the determinant as the product of eigenvalues.
We may use the asymptotical behavior of
empirical spectral distributions, i.e. convergence to the Mar\v
cenko-Pastur distribution in the Wishart/Laguerre case and to the
generalized McKay distribution in the Jacobi case.
However, this structure
is not
"dynamic": if you change $r$, the whole set of eigenvalues is
changing. When using this approach, we will refer to it as the
spectral method.

 In Section 2 we set the framework.
We begin with the matrix models (Wishart-Laguerre, Uniform Gram and
Jacobi), and proceed with the $\beta$-models and processes of
determinants. 

The main results of this paper are in Section 3: laws of large
numbers and fluctuations,   large deviations
and variational problems. The comparison of results obtained by the two methods (decomposition and spectral) deserves interest and is the topic of
Section 4. Some extensions to other models are given in Subsection 4.4.

Sections 5, 6 and 7 are devoted to the proofs.
In Appendix 1 we gather some details on
the Binet formula on the Gamma function which are of constant use in this paper, and Appendix 2 gives
identification of the McKay distribution.

\section{Notation and known facts}
In this long section, we present our different models whose common feature is to introduce 
 processes of random determinants with independent multiplicative factors.  The distribution of these factors  are recorded in Proposition \ref{celebrated} for real matrix models, and settled in formulae (\ref{pLe}), (\ref{pGe}) and (\ref{pJe}) for the (other) $\beta$-models. 

Throughout, $|A|$ stands for $\det A$, and $I_n$  for the $n\times n$ identity matrix. 
If $X$, $Y$ are real random variables and $\mu$ a distribution on $\mathbb R$, we write
\[X \el Y \ \ \ (\hbox{resp.}\ X\el \mu)\] if $X$ and $Y$ have the same distribution (resp. if the distribution of $X$ is $\mu$).

\subsection{Real matrix models and Bartlett-type theorems}
\label{btype}
 In the basic model, we consider  independent random
vectors $b_i, i \geq  1$ with the same distribution $\nu_n$ in
$\mathbb {R}^n$. The most important example is  the Gaussian one with
 $\nu_n = {\mathcal N}(0 ,
I_n)$. If $B = [b_1, \dots, b_r]$, all the entries of $B$ are
independent ${\mathcal N}(0 , 1)$
  and  the distribution of $W=B'B$ is denoted by $W_r(n, \mathbb{R})$ and called  the Wishart ensemble.
  For $r \leq n$,
  its density on the space ${\mathcal S}_r$ of symmetric positive matrix is
\be \frac{1}{2^{rn/2}\Gamma_r(n/2)} |W|^{(n-r-1)/2}
\exp\Big(-\frac{1}{2} \hbox{tr}\!\ W\Big) \ee where $\Gamma_r$ is
the multivariate Gamma function
$$\Gamma_r (\alpha) = \pi^{r(r-1)/4}\Gamma(\alpha)\Gamma\Big(\alpha-\frac{1}{2}\Big)\cdots \Gamma\Big(\alpha-\frac{r-1}{2}\Big)\,,$$

It is the matrix variate extension of the Gamma distribution. Recall
that for $a , c >0$, the $\hbox{Gamma}(a,c)$ distribution has
density
$$\frac{c^a }{\Gamma(a)}\!\ x^{a-1} e^{- cx}\ \ (x > 0)\,.$$
 For $r> n$, the matrix is singular.

Motivated by Hadamard inequality (\ref{Had}), we may choose $\nu_n $ to be
the uniform distribution on the unit sphere 
$\mathbb{S}_{\mathbb R}^{n}$. The corresponding
ensemble for $B$ is called Uniform Spherical Ensemble by 
\citet{Donoho1}. The matrix ensemble for $B'B$ is called the Gram
ensemble by  \citet*{Spince1}, since
$B'B$ is the Gram matrix built from the $b_i$'s. To stress on the
distribution, we call it Uniform Gram ensemble. The diagonal entries
are one and for $r \leq n$, the joint density of the non-diagonal
entries $(r_{ij}, \ 1 \leq i<j\leq r)$ of the matrix $G = B'B$ is
\ben \label{denscorr}
\frac{\left[\Gamma\left(n/2\right)\right]^r}{\Gamma_r\left(n/2\right)}
\ |G|^{(n-r-1)/2} \ \ \ (-1 < r_{ij} < 1) \een (see 
\citet*{GuptaNagar} Theorem 3.3.24 p.107,  \citet{Mathaibook}
Example 1.25
 and \citet{Mathairebook} p.58).

Let us now introduce Jacobi ensembles. For $n_1, n_2 \geq 1$ and $r
\leq n :=n_1 + n_2$, we can decompose every $(n_1+n_2)\times r$ matrix $M$
 in two blocks  \[M= \left(\!\!
\begin{array}{c}
M_1\\
M_2 \end{array}\!\!\right)\] with $M_1$ of type $n_1\times r$ and
$M_2$ of type $n_2 \times r$.
 If the entries of $M$ are
independent ${\mathcal
N}(0,1)$,
 then $W_1 := M_1'M_1$ and $W_2 := M_2'M_2$
are independent Wishart matrices of distribution $W_r(n_1,\mathbb {R})$ and
$W_r(n_2, \mathbb {R})$, respectively. It is well known that $W_1  + W_2$ is
$W_r(n_1 + n_2, \mathbb {R})$ distributed and  a.s. invertible. Let us denote
by $(W_1 + W_2)^{1/2}$   the symmetric positive square root of $(W_1
+ W_2)$.
The  $r\times r$ matrix
\be{\mathcal  X} := (W_1 + W_2)^{-1/2} W_1 (W_1 + W_2)^{-1/2}\ee 
 has a distribution  denoted by $J_r (n_1,
n_2 , \mathbb {R})$ and called the  Jacobi ensemble.

If $T$ is upper triangular with positive diagonal entries and $W_1
+ W_2 = T'T$ (Cholesky decomposition) then \be {\mathcal  Z} =
(T')^{-1} W_1 T^{-1}\ee is also $J_r (n_1, n_2 , \mathbb {R})$ distributed,
(see \citet{Olkin}, \citet{Muir} p.108).

Another occurrence of the Jacobi ensemble is interesting (see
 \citet{Doum}, \citet{Collins}). If $M$ is as above,
its singular value decomposition is \[M=UDV\ \ , \ \ D= \left(\!\!
\begin{array}{c}
\Delta\\
0 \end{array}\!\!\right)
 \] with $D$ of type $n\times r$, with
 $\Delta$ diagonal with nonnegative entries, with
 $U \in \mathcal{O}(n)$ and $V\in \mathcal{O}(r)$ (the orthogonal groups).
 Although $U$ and $V$ are not uniquely determined,
one can choose them according to the  Haar distribution on their respective group and such that $U,V,\Delta$
are independent. Then $M'M= V'\Delta^2V$ and
$$(W_1 + W_2)^{1/2} = (M'M)^{1/2} = V'\Delta V\,.$$
Let $Y_r = U^{[n_1, r]}$ be the $n_1\times r$ upper-left corner of $U$. Since $M_1 = Y_r\Delta V$ we have
$$M_1'M_1 = V'\Delta Y_r'Y_r\Delta V = (MM^*)^{1/2} (V'Y_r'Y_r V)(MM^*)^{1/2}$$
and then ${\mathcal  X}= (Y_r V)'(Y_r V) \el Y_r' Y_r$. In other
words, \be {\mathcal Y} := \big(U^{[n_1, r]}\big)' U^{[n_1, r]} \ee
is also $J_r (n_1, n_2 , \mathbb {R})$ distributed.

If $r\leq \min (n_1, n_2)$, the distribution $J_r (n_1, n_2 , \mathbb {R})$
has a density  on ${\mathcal S}_r$
    which is
\ben \label{betarr} \frac{1}{\beta_r\left(\frac{n_1}{2},
\frac{n_2}{2}\right)} |{\mathcal  Z}|^{\frac{n_1 -r -1}{2}} \ |I_r -
{\mathcal Z}|^{\frac{n_2 -r -1}{2}} \ \mathbf{1}_{0< {\mathcal  Z} <
I_r}\,, \een 
 where
\[\beta_r\left(a, b\right) = \frac{\Gamma_r \left(a\right)\Gamma_r
\left(b\right)} {\Gamma_r \left(a+b\right)}\,,\] (see for example
 \citet{Muir} Theorem 3.3.1). It is the matrix variate
extension of the beta distribution. Recall that for $a > 0, \ b >0$,
the Beta$(a ,b)$ distribution has density \ben \label{recallbeta}
\frac{\Gamma(a + b)}{\Gamma(a)\Gamma(b)} \!\ x^{a -1} (1-x)^{b -1} \
\ (x > 0)\,. \een

Until now, we had $r$ fixed. Our purpose is now to consider all
values of $r$ simultaneously to give a "sample path" study of
determinants.

For  an $n\times n$ matrix $B = [b_1, \dots, b_n]$, we have for $r \leq n$ \be(B'B)^{[r]} =
\Big(B^{[n,r]}\Big)' B^{[n,r]}\,, \ee and for every $j\leq n$, the quantity
\ben\label{rhorhoprem} \rho_{j,n} := \frac{|(B'B)^{[j]}|}{|(B'B)^{[j-1]}|} \een is a
measurable function of $(b_1, \dots, b_j)$ and \ben\label{mesW} |(B'B)^{[r]}| =
\prod_{j=1}^r \rho_{j,n}\,. \een The same occurs  with $\wt b_i :=b_i / \Vert b_i\Vert$ instead of $b_i$ $(i=1,
\dots,n)$ and $\wt B := [\wt b_1 , \dots, \wt b_n]$ instead of $B$.
Let us note that $\wt \rho_{1,n} =1$ and \ben \label{rhorhowt} \wt
\rho_{j,n} =\frac{|\wt W^{[j]}|}{|\wt W^{[j-1]}|} =
\frac{|W^{[j]}|}{|W^{[j-1]}| W_{jj}} = \frac{\rho_{j,n}}{\Vert
b_j\Vert^2} , \ j=2, \dots, n\,,\een
so that
\ben
\label{mesG}
|(\widetilde B'\widetilde B)^{[r]}| =
\prod_{j=1}^r \widetilde\rho_{j,n}\,.
\een

The Wishart
 case and the Uniform Gram case corresponds to (\ref{mesW}) and (\ref{mesG}) respectively, for  $r= 1, \dots,n$.

In the Jacobi case, $r\in \{1,\dots, n_1\}$. If $M = [b_1, \dots,
b_{n_1}]$
,  and if $T$, $W_1$,
${\mathcal Z}$ are defined as above with $n_1$ instead of $r$, then
\be {\mathcal Z}^{[r]} = \Big(\Big(T^{[r]}\Big)'\Big)^{-1} W_1^{[r]}
\Big(T^{[r]}\Big)^{-1}\,. \ee For every $j$, the quantity \be
\rho_{j, n_1, n_2}^{\mathcal Z} := \frac{|{\mathcal
Z}^{[j]}|}{|{\mathcal Z}^{[j-1]}|} \ee is a measurable function of
$(b_1, \dots, b_j)$ and \be |{\mathcal Z}^{[r]}| = \prod_{j=1}^r
\rho_{j,n_1, n_2}^{\mathcal Z}\,. \ee It can be noticed  that \be
\rho_{j,n_1, n_2}^{\mathcal Z} = \frac{|W_1^{[j]}|}{|W_1^{[j]} +
W_2^{[j]}|} \times \frac{|W_1^{[j-1]} +
W_2^{[j-1]}|}{|W_1^{[j-1]}|}\,. \ee

Besides, the construction with the symmetric square root is
different. If \be \rho_{j, n_1, n_2}^{\mathcal X} :=
\frac{|{\mathcal X}^{[j]}|}{|{\mathcal X}^{[j-1]}|} \ee we have
\[{\mathcal X}^{[r]} \not= \Big(W_1^{[r]}+ W_2^{[r]}
\Big)^{-1/2}W_1^{[r]}\Big(W_1^{[r]}+ W_2^{[r]} \Big)^{-1/2}\,.\]
(Take $n_1=n_2 = 2$, $W_1= I_2$, $W_2 = \left(\!\!
\begin{array}{cc}
1 & s\\
s& 1 \end{array}\!\!\right)$ and $r=1$ then
$$\Big(W_1^{[1]}+ W_2^{[1]} \Big)^{-1/2}W_1^{[r]}\Big(W_1^{[r]}+
W_2^{[r]} \Big)^{-1/2} = 2/(4-s^2)\,,$$ and ${\mathcal X}^{[1]} =
1/2$). Moreover we
 cannot say that $\rho_{j, n_1, n_2}^{\mathcal X}$ is measurable
with respect to $b_1, \dots, b_j$ .

Let us consider  the construction from contraction of Haar matrices.
Since \[\Big(\Big(U^{[n_1]}\Big)'U^{[n_1]}\Big)^{[n_1, r]} =
\Big(U^{[n_1, r]}\Big)'U^{[n_1, r]}\,,\] we see that the quantity
\be \rho^{\mathcal Y}_{n_1,n_2,j} := \frac{\left|{\mathcal
Y}^{[j]}\right|}{\left|{\mathcal Y}^{[j-1]}\right|} \ee depends only
on the $j$ first columns of the matrix $U$, and \be |{\mathcal
Y}^{[r]}| = \prod_{j=1}^r \rho^{\mathcal Y}_{n_1,n_2,j}\,.\ee

It is possible to introduce a probability space on which all Uniform
Gram and Wishart matrices are defined for all values of $n$
simultaneously.
   It is enough to consider the infinite product space
   generated by a double infinite sequence of  independent ${\mathcal N}(0,1)$ variables
    $\{b_{i,j}\}_{i,j =1}^\infty$, and for every $n$ to perform the
    above constructions with $b_i = (b_{1i}, \dots, b_{n,i})'$.
    To embed the Jacobi matrices in this framework, we have to
    restrict ourselves to the ${\mathcal X}$-type and ${\mathcal
    Z}$-type ones; however, only the ${\mathcal Z}$ one gives a natural
    meaning to the dynamic study.

The starting point of our study of random determinants is the
following proposition which gathers known results about the factors
entering in the above decompositions.

\begin{prop}
\label{celebrated}
\begin{itemize}
\item[1)] (Bartlett)
The random variables $ \rho_{j,n} , \ j= 1, \dots , n$ are
independent and
\be\rho_{j,n} \el \hbox{Gamma}\Big(\frac{n-(j-1)}{2}, \frac{1}{2}\Big)
\,,\ee
 where $\el$ stands for equality in distribution.
 \item[2)]
The random variables $\wt \rho_{j,n}, \ j=2, \dots, n$ are
independent and \be \wt \rho_{j, n} \el
\hbox{Beta}\Big(\frac{n-j+1}{2}, \frac{j-1}{2}\Big)\,. \ee
 \item[3)] For $J = {\mathcal X}$ (resp. ${\mathcal Y}$, ${\mathcal Z}$), the random variables $\rho^J_{j, n_1, n_2}$ , $j=1, \dots ,
 n_1$
 are independent and
\be 
 \rho^J_{j, n_1, n_2} \el
\hbox{Beta}\Big(\frac{n_1-j+1}{2}, \frac{n_2}{2}\Big)\,. \ee
\end{itemize}
\end{prop}

The first claim is known as the celebrated Bartlett decomposition
(stated with $\chi^2$ distributions) \citep{Bart}. It is
quoted in many books and articles in particular 
\citet{Anderson} pp.170-172, \citet{Muir} Theorem 3.2.14 p.99,
 \citet{Kshi2},  \citet*{GuptaNagar} Theorem
3.3.4 p.91 and ex. 3.8 p.127. The second claim may be
found in \citet{Anderson} Theorem
9.3.3. In the third claim, we first note that it is enough to get
the proof for ${\mathcal Z}$ since the three random matrices have
the same distribution. It is
a consequence of a result quoted in \citet{Anderson}, due to
Kshirsagar, is proved in \citet{Muir} Theorem 3.3.1 p.110
under the assumption $r \leq n_1, n_2$ and in \citet{Rao} p.541
under the only assumption $r\leq n_1$. Actually (see \citet{Muir} ex. 3.24 and  
\citet{Anderson} Theorem 8.4.1), some proofs use probabilistic arguments (as
\citet{Rao} and \citet{Anderson}), Jacobian arguments (as in
\citet*{GuptaNagar} Theorem 5.3.24 p.181), or Mellin transform
arguments (as in \citet{Mathai} Theorem 2).

\subsection{Distribution of eigenvalues and $\beta$-models}
In the study of sttionary processes, random matrices of the Wishart type with complex entries play an important role (\citet{Goodman}).
Less often, quaternionic entries are considered (see \citet{HSS}).
 We do not give details on the complex and quaternionic cases but jump to a general framework. Popularized by physicists, the modern point of
view consists in introducing a parameter $\beta$ taking value $1$
when real, $2$ when complex, and $4$ when quaternionic, this
parameter playing the role of an inverse temperature. The above
constructions can be extended when replacing the transpose ($\mathbb{R}$
{\it case}) by the adjoint ($\mathbb {C}$ {\it case}) or the dual ($\mathbb {H}$
{\it case}). Many of the above results are then true when replacing
in displays the factor $1/2$ by the factor $\beta/2$.

 Actually there are two ways to reach the law of
determinants :

a)  directly  from the distribution of matrices, and using the decomposition method quoted above,

b) from the joint distribution of eigenvalues.

The second way has been used to define the so-called $\beta$-models,
see for instance  \citet{Forrester} Chap.2. The idea of
extending the range of  $\beta$ to $(0, \infty)$ is quite natural. As mentioned in Section 1,
they correspond also to models of
tri-diagonal random matrices  (see 
\citet*{EdelmanDumitriu}, \citet*{Killip1}, \citet*{edelman:bjm}).

Because of the connections with orthogonal polynomials in the
complex case, the extended families are called $\beta$-Laguerre
ensemble (or just Laguerre ensemble) instead of Wishart ensemble and
$\beta$-Jacobi ensemble (or just Jacobi ensemble) instead of MANOVA
or Beta ensemble.

Throughout, we use the symbol $\beta'$ for $\beta/2$ to simplify
displays.

\subsubsection{Laguerre}
When $\beta =1,2,4$ the joint probability density of the eigenvalues
$\lambda_j , j = 1, \dots , r$ of $W$ on the orthant  $\lambda_j > 0
 , \ j=1, \dots , r$ is
 \ben
 \label{jdL}
 \frac{1}{Z^{L,\beta}_{r}(n)}\!\ \prod_{j=1}^r\left(
\lambda_j^{\beta'(n-r + 1) -1} \e^{-\beta' \lambda_j}\right)\
\prod_{1 \leq j < k \leq r} |\lambda_k - \lambda_j|^{2\beta'}\ \ , \
\een and the normalizing constant is \be Z^{L,\beta}_{r}(n)=
\Big(\frac{1}{\beta'}\Big)^{\beta' rn} \prod_{j=1}^r
\frac{\Gamma\left(1 + \beta'
j\right)\Gamma\left(\beta'(n-j+1)\right)}{\Gamma\left(1+
\beta'\right)}\,.\ee This is the Selberg integral (see for instance
 \cite{EdelmanRao}, formula 4.6 and references
therein).

 When $\beta > 0$ is not $1,2,4$, we consider the density (\ref{jdL}) on $(0,\infty)^r$. We also
denote the product $\prod_{j=1}^r \lambda_j$ by $|W|$. This gives
the Mellin transform 
\be 
 \mathbb{E} |W|^{\beta' s} =
\frac{Z^{L,\beta}_{r}(n+s)}{Z^{L,\beta}_{r}(n)} =
\Big(\frac{1}{\beta'}\Big)^{\beta' rs} \ \prod_{k=1}^r
\frac{\Gamma\left(\beta'(n-k+1+s)\right)}{\Gamma\left(\beta'(n-k+1)\right)}\,.
\ee Remembering that if $X \el \hbox{Gamma} (a, 1/2)$ then
\[\mathbb{E} X^\mu = 2^\mu \frac{\Gamma(\mu + a)}{\Gamma(a)} \ \ (\mu
> -a)\,,\] we deduce the following proposition from the uniqueness of Mellin transform.
\begin{prop}
\label{pLe}We have \be
|W| \el \prod_{j=1}^r \rho^{L, \beta}_{j,n}\,,  \ee where the variables
$\rho^{L, \beta}_{j,n}, j = 1, \dots, r$ are independent and \ben
\label{lawrhoL} \rho^{L, \beta}_{j,n} \el
\hbox{Gamma}\left(\beta'(n-j+1) , 1/2\right)\,.
 \een
 \end{prop}
 We stress that our point of view
is not compatible with the construction by 
(\cite{EdelmanDumitriu}) of matrix models for the (general)
$\beta$-Laguerre ensemble. Actually, they define a random $r\times r$
matrix $B^{(r)}$ where only
 diagonal and subdiagonal terms are  nonzero, independent and satisfy (for $n$ fixed):
\be B^{(r)}_{ii} &\el& \sqrt{Gamma\left(\beta'(n-i+1),
1/2\right)}\ \ \ (1 \leq i\leq r)\,,\\B^{(r)}_{i, i-1} &\el&
\sqrt{Gamma\left(\beta'(r-i+1), 1/2\right)}\ \ \ (2 \leq
i\leq r)\,. \ee They prove that the distribution of eigenvalues of
$B^{(r)} \Big(B^{(r)}\Big)'$ is precisely (\ref{jdL}). Of course we
recover the determinant as a product of elements with the good
distribution, but the problem is that we cannot consider all $r$
simultaneously in their framework, since
\[\Big(B^{(r)}\left(B^{(r)}\right)'\Big)^{[r-1]} \not= B^{(r-1)}\Big(B^{(r-1)}\Big)'\,.\]

\subsubsection{Uniform Gram}
It is useful in the study of correlations. A correlation matrix is a
positive definite matrix with diagonal entries equal to one. Here,
there is no explicit expression for the law of eigenvalues.
However, the expression \be  \frac{1}{Z^{G,\beta}_{r}(n)} \!\
|G|^{\beta' (n - r+1)-1} \ee with \be Z^{G,\beta}_{r}(n)
=\pi^{\beta' r(r-1)}\!\ \prod_{j=1}^r \frac{\Gamma\left(\beta'
(n-j+1)\right)}{\Gamma\left(\beta' n\right)} \ee
 is a density on the space of symmetric (resp. Hermitian, resp.
self-dual) positive matrices with diagonal entries equal to one, and
it fits with the distribution of correlation matrix in the real (see
(\ref{denscorr})), complex and quaternion case, for the appropriate values of $\beta$. This yields
(\cite{GuptaNagar} ex. 3.26 p.130)  the Mellin transform \be \mathbb{E}
|G|^{\beta' s}  = \frac{Z^{G,\beta}_{r}(n+s)}{Z^{G,\beta}_{r}(n)} =
\prod_{j=1}^r
\frac{\Gamma\left(\beta'(n-j+1+s)\right)\Gamma\left(\beta'
n\right)}{\Gamma\left(\beta'(n-j+1)\right)\Gamma\left(\beta' (n
+s)\right)}
 \ee
From (\ref{recallbeta}), it is clear that if $X\el \hbox{Beta}(a,b)$
then \ben \label{Mellinbeta} \mathbb{E} X^\mu = \frac{\Gamma(a + \mu)
\Gamma(a+b)}{\Gamma(a) \Gamma(a+b+\mu)} \ \ \ (\mu > -a)\,. \een Again
the uniqueness of the Mellin transform leads to the proposition.
\begin{prop}
\label{pGe} We have
 \be |G|
\el \prod_{j=2}^r \rho^{G, \beta}_{j,n} \ee where the variables
$\rho^{G, \beta}_{j,n}, j = 2, \dots, r$ are independent and \ben
\label{lawrhoG} \rho^{G, \beta}_{j,n} \el
\hbox{Beta}\left(\beta'(n-j+1) , \beta'(j-1)\right)\,. \een
\end{prop}

\subsubsection{Jacobi}

If ${\mathcal Z}$ is distributed as in (\ref{betarr}), the joint
density of eigenvalues on the set $(0 < \lambda_j < 1 \, \ j=1,
\dots, r)$ is given by \be
 \frac{1}{Z_r\left(n_1,
n_2\right)}\prod_{i=1}^r \lambda_i^{\frac{n_1 - r -1}{2}} \left(1 -
\lambda_i\right)^{\frac{n_2 - r -1}{2}}\prod_{1 \leq i <j \leq r}
|\lambda_j - \lambda_i| \,, \ee
where $Z_r\left(n_1, n_2\right)$ is a normalizing constant
 (see for example \cite{Muir} Theorem 3.3.4).

For $n_2 < r < n_1$, the matrix $W_2$ is singular, and  the Jacobi
matrix $I - {\mathcal Z}$
 has $1$ as an eigenvalue with multiplicity $r-n_2$. The
distribution of ${\mathcal Z}$ has no density. Nevertheless we may
study its determinant. Indeed, the matrix $I - {\mathcal Z}$
 has $0$ as an  eigenvalue of
multiplicity $r-n_2$. Actually  the density of the law of the
non-zero eigenvalues of this matrix is known (see 
\cite{Sri1} and  \cite{DIAZ}), so
that the non-one eigenvalues of ${\mathcal Z}$ have the joint
density \be 
 \frac{1}{\wt Z_r\left(n_1,
n_2\right)}\prod_{i= 1}^{n_2} \lambda_i^{\frac{n_1 - r -1}{2}}
\left(1 - \lambda_i\right)^{\frac{r -n_2  -1}{2}}\prod_{  1 \leq i
<j \leq n_2} |\lambda_j - \lambda_i|\,, \ee
where the normalizing
constant is $\wt Z_r\left(n_1, n_2\right) = Z_{n_2} (n_1 + n_2 -r,
r)\,.$

We now consider matrices with elements in $\mathbb {X} = \C$ or $\H$. When $r
\leq \min(n_1, n_2)$, the 
distribution
 of ${\mathcal Z}$ has a density proportional to
\be
 |{\mathcal  Z}|^{\beta'(n_1 -r +1)-1}\ |I - {\mathcal  Z}|^{\beta'(n_2 -r +1)-1} \ \mathbf{1}_{0< {\mathcal  Z} < I} \ \,.
\ee where $\beta'=1$ or $2$. The distribution  of the 
eigenvalues of ${\mathcal Z}$ has the density (on $[0,1]^r$) :
\ben
\label{vpjbetab}
f^\beta_{r,n_1, n_2}(\lambda_1, \dots , \lambda_r)= \ \ \ \ \ \ \ \ \  \ \ \ \ \ \ \ \ \
 \ \ \ \ \ \ \ \ \ \ \ \ \ \ \ \ \ \  \ \ \ \ \ \ \ \ \  \ \ \ \ \ \ \ \ \\\
\nonumber  \frac{1}{Z_r^{(J, \beta)}\left(n_1,
n_2\right)}\prod_{i=1}^r \lambda_i^{\beta'(n_1 -r +1)-1} \left(1 -
\lambda_i\right)^{\beta'(n_2 -r +1)-1} \prod_{1 \leq i <j \leq r}
|\lambda_j - \lambda_i|^{2\beta'}\,,\een
 where
\ben \label{Selberg} Z_r^{(J, \beta)} (n_1, n_2) = \prod_{j=1}^{r}
\frac{ \Gamma\left(1+ \beta' j\right)\Gamma\left(\beta'(n_1
+j-r)\right) \Gamma\left(\beta'(n_2 +j-r)\right) } {\Gamma\left(1+
\beta'\right)\Gamma\left(\beta'(n_1+n_2 +j-r)\right)}\,, \een is the
value of the Selberg integral
 (see  \cite{HiaiP} p.118 and
also  \cite{EdelmanRao} p.19 and references
therein).

In the "singular" case $(n_2 \leq r \leq n_1)$, the density of the
non-one eigenvalues is $f^\beta_{n_2, n_1 +n_2 -r, r}(\lambda_1,
\dots, \lambda_{n_2})$.

We consider an extension of the above models. For every $\beta > 0$,
we define a family of distribution densities ${\bf f}_{r,n_1,
n_2}^\beta$  on $[0,1]^{\min(n_2 , r)}$ :

\begin{equation}
\label{defbff}
 {\bf f}_{r,n_1, n_2}^\beta =
     \begin{cases}
       f_{r,n_1, n_2}^\beta   & \text{if $r \leq \min(n_1 , n_2)$} \\
       f_{n_2, n_1 +n_2 -r, r}^\beta & \text{if $n_2 \leq r \leq n_1$}\,.
     \end{cases}
\end{equation}

 We set by convention $$|{\mathcal  Z}_{n_1, n_2, r}| = \prod_{i=1}^{\min(n_2 , r)} \lambda_i$$ in all cases,  and we call it the determinant, even if we do not define any matrix.

For $r \leq n_1, n_2$, using (\ref{vpjbetab}) and (\ref{Selberg}) we
obtain \ben \nonumber \mathbb{E} \left(| {\mathcal Z}_{n_1, n_2, r}|^{\beta'
s}\right) = \frac{Z_r^{J,\beta} \left(n_1 + s, n_2
\right)}{Z_r^{J,\beta}\left(n_1, n_2\right)}\ \ \ \ \ \ \ \ \ \ \ \
\ \ \ \ \ \ \ \ \ \ \ \ \ \ \ \ \ \ \ \ \ \ \ \
\\\label{melreg} = \prod_{j=1}^r \frac{ \Gamma\left(\beta'(n_1 +n_2
+j-r)\right)\Gamma\left(\beta'(n_1 +j-r +s)\right) }
{\Gamma\left(\beta'(n_1 +j-r)\right)\Gamma\left(\beta'(n_1+n_2
+j-r+s)\right)}\,. \een

If $n_2 < r \leq n_1$ we start directly from (\ref{defbff}) and
(\ref{Selberg}) we have \ben \nonumber  \mathbb{E} \left(| {\mathcal
Z}_{n_1, n_2,r}|^{\beta' s}\right) &=& \frac{Z_{n_2}^{J,\beta}
\left(n_1 +n_2 -r + s, r \right)}{Z_{n_2}^{J, \beta}\left(n_1+n_2
-r, r \right)}\ \ \ \ \ \ \ \ \ \ \ \ \ \ \ \ \ \ \ \ \ \ \ \ \ \ \
\ \ \ \ \ \ \ \ \
\\ \label{melregs}&=& \prod_{j=1}^{n_2} \frac{ \Gamma\left(\beta'(n_1
+j)\right)\Gamma\left(\beta'(n_1 +j-r+s)\right) }
{\Gamma\left(\beta'(n_1 +j-r)\right)\Gamma\left(\beta'(n_1
+j+s)\right)}\,. \een Multiplying up and down by $\prod_{k= n_2
+1}^r \Gamma\left(\beta'(n_1 +k-r)\right) \Gamma\left(\beta'(n_1
+k-r+s)\right)$ we get again the right hand side of (\ref{melreg}).
Going back to
 (\ref{Mellinbeta}), we have the proposition
 \begin{prop}
 \label{pJe}
For $r \leq n_1$ :
$$|{\mathcal  Z}_{n_1, n_2, r}|\el \prod_{j=1}^{r} \rho^{\beta,J}_{j,n_1, n_2}\,,$$
where $\rho^{\beta,J}_{j,n_1, n_2} , \ j= 1, \dots , r$ are
independent and \ben \label{loih} \rho^{\beta,J}_{j,n_1, n_2} \el
\hbox{Beta}\left(\beta'(n_1-j+1), \beta' n_2\right)\,. \een
\end{prop}

 \subsection{Processes}
In the three ensembles defined above, we have met arrays of independent variables
with remarkable distributions. In Section \ref{btype}, we have
discussed the interest of studying all values of $r$ simultaneously
in the matrix cases ($\beta = 1,2,4$). Since the structure remains
the same  in the $\beta$-models, it is meaningful to consider the
processes (indexed by $r$) of partial sums. A now classical
asymptotic regime is $n,r \rightarrow \infty$ with fixed ratio in
the Laguerre and Uniform Gram case, and  $n_1, n_2, r \rightarrow
\infty$ with fixed ratios in the Jacobi case. It means that we
consider the asymptotic behavior determinants in a dynamic (or path
wise) way.

For the Laguerre model, we define \ben\label{sumL} \log \Delta^{L,
\beta}_{n,p} := \sum_{k=1}^p \log \frac{\rho^{L, \beta}_{k,n}}{\beta
n} \ \ \  \ (p \leq n)\een and the process \ben \label{pL}
\Delta^{L, \beta}_{n}(t) := \Delta^{L, \beta}_{n, \lfloor
nt\rfloor}, \ \ t \in [0, 1]\,. \een
 For the Uniform Gram
model, we define \ben\label{sumG}  \log \Delta^{G, \beta}_{n,p} :=
\sum_{k=1}^p \log \rho^{G, \beta}_{k,n}\ \ \  \ (p \leq n)\een and the
process \ben \label{pG} \Delta^{G, \beta}_{n}(t) := \Delta^{G,
\beta}_{n, \lfloor nt\rfloor}, \ \ t \in [0, 1]\,.\een

For the Jacobi model, we fix $\tau_1$ and $\tau_2
> 0$, set $n_1 = \lfloor n\tau_1\rfloor, n_2 = \lfloor
n\tau_2\rfloor$, and define \ben\label{sumJ} \log \Delta^{J,
\beta}_{n, p} = \sum_{k=1}^p \log \rho^{J, \beta}_{k,n_1, n_2}\ \ \  \ (p \leq n_1)\een and the process \ben \label{pJ} \Delta^{J, \beta}_n
(t) = \Delta^{J, \beta}_{n, \lfloor nt\rfloor}, \ \ t \in
[0,\tau_1]\,.\een

There are some connections between the above processes. For
instance,  in the real matrix model ($\beta=1$) we saw in
(\ref{rhorhowt}) that
$$\rho^{L, 1}_{j,n} = \rho^{G, 1}_{j,n}\ \Vert b_j\Vert^2\,, $$
that the two random variables in the right hand side are independent
and $\Vert b_j\Vert^2 \el \hbox{Gamma}(n/2 , 1/2)$.

To see these connections in the general case, we use  the so-called
"beta-gamma" algebra that will be really helpful in the sequel. Details
can be found in \cite{ChauYor} pp.93-94. In the following relation,
$\gamma(a)$ denotes a random variable with distribution
Gamma$(a,1)$, and $\beta(a,b)$ denotes a random variable with
distribution Beta$(a,b)$. The relation is \ben \label{gamma0}
\big(\gamma(a) , \gamma(b)
\big)&\sur{=}{(d)}&\big(\beta(a,b)\gamma(a+b) , (1 - \beta(a,b))
\gamma(a+b)\big)\,, \een where, on the left hand side the random
variables $\gamma(a)$ and $\gamma(b)$ are independent and on the
right hand side the random variables $\beta(a,b)$ and $\gamma(a+b)$
are independent. It entails in particular \ben \label{gamma1}
\frac{\gamma(a)}{\gamma(a) + \gamma(b)}&\sur{=}{(d)}& \beta(a,b)\,.
\een Let us note that this relation can be extended at the matrix
variate level.

 From the definitions (\ref{sumL}) and (\ref{sumG})
and owing to the equalities in distribution (\ref{lawrhoL}) and
(\ref{lawrhoG}), we have then \ben \label{relproc} \log\Delta_n^{L,
\beta} \el  \log \Delta^{G, \beta}_n + S_n
 \,,
 \een
where  $S_n$  is independent of  $\log
\Delta^{G, \beta}_n$, and specified by \ben \label{defsigmaproc}
S_n(t) = \sum_{k=1}^{\lfloor nt\rfloor} \log \varepsilon_k^{(n)}\ \
, \ t\in [0,1]\een where $\varepsilon_k^{(n)}$, $k=1, \dots , n$ are
independent and satisfy $\varepsilon_k^{(n)} \el
\hbox{Gamma}\left(\beta' n , \beta' n\right)\,.$  In the sequel, we
begin by setting the claims for the Uniform Gram process and then
deduce the corresponding results for the Laguerre process.

 Using the definitions
(\ref{sumL}) and (\ref{sumJ}) and  the equalities in distribution
(\ref{lawrhoL}) and (\ref{loih}), we get, by another application of
(\ref{gamma0})

 \ben
 \label{ThetamoinsTheta}
\log \Delta_{n_1, r}^{L, \beta} \el \log \Delta_{n, r}^{J, \beta} +
\log \Delta_{n_1 + n_2, r}^{L, \beta}- r \log\frac{n_1}{n_1 +
n_2}\,,
 \een
 where this equality holds for all indices $r=1, \dots, n_1$
 simultaneously, and the two processes $\log \Delta_{n}^{J, \beta}$ and $\log \Delta_{n_1 + n_2}^{L, \beta}$ are independent.

It allows to deduce asymptotic results for the Jacobi model from
those of the Laguerre model.

\section{Main results}
In this section, we state first a law of large numbers and
fluctuations for our three models, and then  the corresponding LDP
for processes and marginals.

 Let $D_T = \{ v \in \mathbb{D} ([0,T]) : v(0) =
0\}$ the set of c\`adl\`ag functions on $[0,T]$
 and  $D= \{ v \in \mathbb {D}([0,1))  : v(0) = 0\}$
the set of c\`adl\`ag functions on $[0,T]$ and $[0, 1)$, respectively, starting from $0$.

We use often the following entropy function
\begin{equation}
\label{defJ}
 {\mathcal J}(u) =
     \begin{cases}
         u\log u -u + 1 & \text{if $u > 0$} \\
      1  & \text{if $u=0$}\\
  +\infty   & \text{if $u<0$}
     \end{cases}
\end{equation}
and its primitive: \ben \label{defF} F(t) = \int_0^t {\mathcal
J}(u)\ du = \frac{t^2}{2}\log t -\frac{3t^2}{4} + t, \ \ \ (t \geq
0)\,. \een We  use also the function defined in 
\cite{hiai2}, for $s,t \geq 0$: \ben \label{b1}\nonumber
B(s,t)&:=&  \frac{(1+s)^2}{2}\log (1+s) - \frac{s^2}{2}\log s + \frac{(1+t)^2}{2}\log (1+t)  - \frac{t^2}{2}\log t\\
&-& \frac{(2+s+t)^2}{2}\log (2+s+t) + \frac{(1+s+t)^2}{2}\log
(1+s+t)\,. \een which may also be written as
\be
B(s,t) = F(1+s) -F(s) + F(1+t) -F(t) -F(2+s+t) + F(1+s+t) -
\frac{7}{4}\,. \ee
\subsection{Law of large numbers and fluctuations}
\label{LLNF}
\subsubsection{Uniform Gram ensemble}
Let us define a drift and a diffusion coefficient by \ben
\label{defdriftG} \d^{G, \beta}(t) := \frac{1}{\beta} +
\left(\frac{1}{2} - \frac{1}{\beta} \right) \frac{1}{1-t} \ , \
\sigma^{G, \beta}(t) := \sqrt{\frac{2t}{\beta(1-t)}}\,.
 \een
\begin{theorem}
\label{first2}
\begin{enumerate}
\item As $n \rightarrow \infty$,
\ben
 \label{enough} \lim_n \
\sup_{p \leq n} \left|
 \frac{1}{n}\mathbb{E}\!\  \log \Delta^{G, \beta}_{n , p} + {\mathcal J}\left(1 - \frac{p}{n}\right)\right| =
 0\,.
\een
\item For every $ t \in [0,1)$, as $n \rightarrow \infty$,
\ben  \label{esptG}  \mathbb{E} \!\ \log \Delta^{G, \beta}_{n}(t) +
n{\mathcal J}\left(1 - \frac{\lfloor nt\rfloor}{n}\right)
\rightarrow \int_0^t \d^{G, \beta}(s)\!\ ds
 \een and \ben
 \label{esp1G}
   \ \mathbb{E} \!\  \log \Delta^{G, \beta}_{n}(1) +n +  \left(\frac{1}{\beta} -\frac{1}{2}\right) \log n
   \rightarrow
K^1_\beta
   \,,\een
where \ben \label{defk1beta} K^1_\beta  :=  \frac{1}{2}\log (2\pi)
+\frac{1- \gamma}{\beta} - \int_0^\infty \frac{sf(s)}{e^{\beta s/2} -
1}\ ds\,, \een and $\gamma= -\Gamma'(1)$ is the Euler constant.
\item
For every $t \in [0, 1)$, as $n \rightarrow \infty$,\ben
\label{vartG}
 \hbox{Var} \  \log \Delta^{G, \beta}_n (t) &\rightarrow& \int_0^t \Big(\sigma^{G,\beta} (s)\Big)^2\ ds
 \\
\label{var1G}
   \ \hbox {Var} \  \log \Delta^{G, \beta}_{n}(1) - \frac{2}{\beta} \log n  &\rightarrow&  K^2_\beta
   \,,
\een where \ben \label{defk2beta} K^2_\beta :=
\frac{2(\gamma-1)}{\beta} + \int_0^\infty \frac{s(sf(s) +
\frac{1}{2})}{e^{\beta s/2} - 1}\ ds\,. \een
\item As $n \rightarrow \infty$,
\ben \label{theocvps} \lim_n \ \sup_{t\in [0,1]} \left|\frac{ \log
\Delta^{G, \beta}_n (t)}{n} + {\mathcal J}(1-t)\right| = 0\,. \een
 in probability.
\end{enumerate}
\end{theorem}

\medskip

For $\beta = 1$, the formulae (\ref{esp1G}) and (\ref{var1G}) are
due to \cite{abbott} (see their lemmas 4.2 and 4.4),
using a variant of the decomposition method.

\begin{theorem}
\label{DonskerH}
\begin{enumerate}
\item Let for $n \geq 1$
\begin{eqnarray*}
\eta_n^{G,\beta} (t) :=  \log \Delta^{G,\beta}_n (t)
 + n{\mathcal J}\left(1 - \frac{\lfloor nt\rfloor}{n}\right)
 \ \ , \ \ t\in [0,1)\,.
 \end{eqnarray*}
Then as $n \rightarrow \infty$ \ben \label{etantG}
\Big(\eta_n^{G,\beta} (t); \ t \in [0,1)\Big) \Rightarrow \Big(
X_t^{G, \beta} ; \ t \in [0,1)\Big)\,, \een where $X^{G,\beta}$ is
the (Gaussian) diffusion solution of the stochastic differential
equation : \ben \label{sdeG} dX^{G,\beta}_t = \d^{G,
\beta}(t)
\!\ dt + \sigma^{G, \beta}(t)
\ d{\bf B}_t\,, \een with $X^{G,\beta}_0 = 0$, ${\bf B}$ is a
standard Brownian motion and $\Rightarrow$ stands for the weak
convergence of distributions in $D$ endowed with
the Skorokhod topology.

\item Let
\begin{eqnarray*}
\widehat\eta_n^{G,\beta}  = {\frac{\log \Delta^{G,\beta}_{n}(1) +n +
\left(\frac{1}{\beta}-\frac{1}{2}\right) \log
n}{\sqrt{\frac{2}{\beta}\log n}}}\ \,.
\end{eqnarray*}
Then as $n \rightarrow \infty$, $\widehat\eta_n^{G,\beta}
\Rightarrow N$ where $N$ is ${\mathcal N}(0,1)$ and independent of
$\bf B$,
(and $\Rightarrow$ stands for the weak convergence of distribution
in $\mathbb{R}$).
\end{enumerate}
\end{theorem}

\subsubsection{Laguerre ensemble}

Let us define a drift and a diffusion coefficient by \ben
\label{defdriftL} \d^{L, \beta}(t) := \left(\frac{1}{2}
-\frac{1}{\beta}  \right) \frac{1}{1-t} \ \ , \ \ \sigma^{L,
\beta}(t) := \sqrt{\frac{2}{\beta(1-t)}}\,.
 \een
\begin{theorem}
\label{ODD}
\begin{enumerate}
\item As $n \rightarrow \infty$,
\ben \label{unifw} \lim_n \sup_{p\leq n} \left|\frac{1}{n}\!\  \E
\log \Delta^{L, \beta}_{n,p}  + {\mathcal J}\left(1 -
\frac{p}{n}\right)\right| = 0 \een
\item For every $ t \in [0,1)$, as $n \rightarrow \infty$,
\ben
\label{esptW}
 \mathbb{E}\, \log \Delta^{L, \beta}_{n}(t)
+ n{\mathcal J}\left(1 - \frac{\lfloor nt\rfloor}{n}\right)
\rightarrow \int_0^t \d^{L, \beta}_L(s)\ ds \,, \een and \ben
\label{esp1W}
 \mathbb{E}\,
 \log \Delta^{L, \beta}_{n}(1)
  +n  + \left(\frac{1}{\beta}-\frac{1}{2}\right) \log n \rightarrow K_\beta^1
  \,,
\een
\item
For every $t \in [0, 1)$, as $n \rightarrow \infty$, \ben
\label{vartW} \hbox{Var}\, \log \Delta^{L, \beta}_{n}(t)
&\rightarrow& \int_0^t \Big(\sigma^{L, \beta}(s)\Big)^2 \ ds
\\
\label{var1W} \hbox{Var}\,  \log \Delta^{L, \beta}_{n}(1)
 - \frac{2}{\beta} \log n   &\rightarrow& K_\beta^2
 \,.
\een
\item As $n \rightarrow \infty$,
\ben \label{cvpW} \sup_{t \in [0, 1]} \left|\frac{1}{n} \log
\Delta^{L, \beta}_n (t) + {\mathcal J}(1- t) \right|\rightarrow 0
\een in probability.
\end{enumerate}
\end{theorem}
\begin{remark}
In the Uniform Gram and Laguerre ensembles, when all the variables
are defined on the same space (i.e. $\beta=1,2,4$), an application
of the Borel-Cantelli lemma leads to almost sure convergence.
\end{remark}

\begin{theorem}
\label{DonskerW}
Let
\begin{eqnarray*}
\eta_n^{L, \beta} (t) &:=& \log \Delta^{L, \beta}_{n}(t)  + n
{\mathcal J}\left(1 - \frac{\lfloor nt\rfloor}{n}\right)
 \ \ , \ \ t\in [0,1)\,,\\
\widehat\eta_n^{L, \beta}  &=& \frac{ \log \Delta^{L, \beta}_{n}(1)
+n  + \left(\frac{1}{\beta}-\frac{1}{2}\right)\log
n}{\sqrt{\frac{2}{\beta}\log n}}\ \,.
\end{eqnarray*}
Then as $n \rightarrow \infty$
\ben
\label{etantW}
\Big(\eta_n^{L, \beta} (t); \ t \in [0,1)\Big) &\Rightarrow& \Big(X_t^{L, \beta}, \ \ t \in [0,1)\Big)\\
\nonumber \widehat\eta_n^{L, \beta} &\Rightarrow& N \een where
$X^{L, \beta}$ is the Gaussian diffusion solution of the stochastic
differential equation: \ben \label{sdeW} d X_t^{L, \beta} =
\d^{L,\beta}(t)\!\ dt + \sigma^{L, \beta}(t)
\ d {\bf B}_t\,, \een with $X^{L, \beta}_0=0$, where ${\bf B}$ is a
standard Brownian motion and $N$ is ${\mathcal N}(0,1)$ and
independent of $\bf B$.
\end{theorem}
The convergence of $\eta_n^{L,1}(t)$, for fixed $t$ and of $\wh\eta_n^{L,1}$ were proved by \cite{jonsson1} Theorem 5.1a.
Recently and independently the convergence of $\widehat\eta_n^{L, 1}$ was proved  in Theorem 4 of \cite{Rempala1}.
\subsubsection{Jacobi ensemble}
In this part we use new auxiliary functions. Let
$${\mathcal E}(x , y ,z) = x \log x - (x+y) \log (x+y)
+(x+y- z) \log (x+y-z)- (x-z) \log (x-z)$$ or using ${\mathcal J}$
defined in (\ref{defJ}) \ben \label{defe} {\mathcal E}(x , y ,z) =
{\mathcal J}(x) - {\mathcal J}(x-z) -{\mathcal J}(x+y) + {\mathcal
J}(x+y -z)\,. \een The partial derivative of ${\mathcal E}$ with
respect to $x$ is : \ben \label{defe'} {\mathcal E}_1(x , y , z) :=
\frac{\partial}{\partial x}\!\ {\mathcal E}(x , y ,z) = \log \frac{x
(x+y-z)}{(x-z)(x+y)}\,. \een
 Let for $0 \leq t < \tau_1$
 \ben
 \label{defsigma}
 \sigma^2(t) := \frac{\partial}{\partial t}\!\ {\mathcal E}_1 (\tau_1, \tau_2,t) = \frac{\tau_2}{(\tau_1 -t)(\tau_1 + \tau_2-t)}\,.
 \een
Again we define drift and diffusion coefficients: \ben
\label{defdriftJ} \d^{J, \beta} (t) = \left(\frac{1}{2}
-\frac{1}{\beta}  \right) \sigma^2(t)\ \ , \ \  \nonumber \sigma^{J,
\beta} (t) = \sqrt{\frac{2}{\beta}}\ \sigma(t)\,.
 \een
\begin{theorem}
\label{thlln}
\begin{enumerate}
\item As $n \rightarrow \infty$,
\ben \label{unif} \sup_{t \in [0, \tau_1]}\left|\frac{1}{n}\mathbb{E}
\log\Delta_{n}^{J, \beta}(t) - {\mathcal E}\left(\tau_1 , \tau_2
,t\right)\right|  \rightarrow 0\,. \een
\item For every $ t \in [0, \tau_1)$, as $n \rightarrow \infty$,
\ben \label{dj} \mathbb{E} \log\Delta_{n}^{J, \beta} (t) - {\mathcal
E}(\lfloor\tau_1 n\rfloor, \lfloor\tau_2 n\rfloor,\lfloor
tn\rfloor)\longrightarrow \int_0^t \d_J (s)\ ds \,, \een
and\footnotemark[1]\footnotetext[1]{where $K^1_\beta$ (resp.
$K^2_\beta$) was defined in (\ref{defk1beta}) (resp.
(\ref{defk2beta}))\label{note}} \ben \nonumber \mathbb{E} \log\Delta_{n}^{J,
\beta} (\tau_1)- {\mathcal E}(\lfloor\tau_1 n\rfloor, \lfloor\tau_2
n\rfloor,\lfloor \tau_1 n\rfloor)+
\left(\frac{1}{\beta}- \frac{1}{2}\right)\log n \longrightarrow\\
\label{djc}  \left(\frac{1}{2}-\frac{1}{\beta}\right) \log
\frac{\tau_1\tau_2}{\tau_1 + \tau_2}+  K^1_\beta\,. \een .
\item
For every $t \in [0, \tau_1)$, as $n \rightarrow \infty$,
\ben
\label{limvar} \hbox{Var}\ \log\Delta_{n}^{J, \beta} (t) \rightarrow
\int_0^t \Big(\sigma^{J, \beta}(s)\Big)^2 \ ds\,,
\een and\footnotemark[1]\ben \label{limvarc} \hbox{Var}\
\log\Delta_{n}^{J, \beta} (\tau_1) -\frac{2}{\beta}\log n
\longrightarrow \frac{2}{\beta} \log \Big(\frac{\tau_1\tau_2}{\tau_1
+ \tau_2}\Big) + K^2_\beta\,. \een
\item As $n \rightarrow \infty$,
\ben \label{esplim} \sup_{t \in [0, \tau_1]} \left|\frac{1}{n}
\log\Delta_{n}^{J, \beta} (t) -{\mathcal E}(\tau_1, \tau_2 ,t)
\right|\rightarrow 0 \een in probability.
\end{enumerate}
\end{theorem}
\begin{remark}
For $\beta = 1, 2,4$, when all  variables are on the same
probability space, the convergence in (4)  may be strengthened to
almost sure convergence.
\end{remark}

\begin{theorem}
\label{DonskerJ} Let for $n \geq 1$
\begin{eqnarray*}
\eta_n^{J, \beta} (t) &:=&  \log\Delta_{n}^{J, \beta} (t) -
{\mathcal E}(\lfloor\tau_1 n\rfloor, \lfloor\tau_2 n\rfloor,\lfloor
tn\rfloor)
  \ \ , \ \ t\in [0,\tau_1)\,,\\
\widehat\eta_n^{J, \beta}  &:=& \frac{\log\Delta_{n}^{J, \beta}
(\tau_1) - n {\mathcal E}(\tau_1, \tau_2, \tau_1)+
\left(\frac{1}{2}-\frac{1}{\beta}\right)\log
n}{\sqrt{\frac{2}{\beta} \log n}}\ \,.
 \end{eqnarray*}
Then as $n \rightarrow \infty$ \ben \label{etantHJ} \Big(\eta_n^{J,
\beta}
 (t); \ t \in [0,\tau_1)\Big) &\Rightarrow& \Big( X_t^J ; \ t \in
 [0,\tau_1)\Big)\,,\\
\nonumber \widehat\eta_n^{J, \beta} &\Rightarrow& N\een where $X^{J,
\beta}$ is the (Gaussian) diffusion solution of the stochastic
differential equation : \ben \label{sdeJ} dX^J_t = \d^{J, \beta}(t)
dt + \sigma^{J, \beta}(t)\ d{\bf B}_t\,, \een with $X^{J, \beta}_0 =
0$, ${\bf B}$ is a standard Brownian motion and $N$ is a standard
normal variable independent of $\bf B$.
\end{theorem}

\subsection{Large deviations}
\label{LD} All along this section, we use the notation of  \cite{DZ}. In particular we write LDP for Large Deviation
Principle. The reader may have some interest in consulting  \cite{gamb} where a similar method is used for a different
model, but here we use  a slightly different topology to be able to
catch the marginals  in $T$.

For $T < 1$, let 
 $\mathrm{M}_T$ be the set of
signed measures on $[0,T]$ and let ${\mathrm M}_<$ be the set of
measures whose
 support is a compact subset of $[0,1)$.

We provide $D$ with the weakened topology $\sigma(D, {\mathrm
M}_<)$. So, $D$ is the projective limit of the family, indexed by
$T< 1$ of topological spaces  $\left(D_T , \sigma(D_T,{\mathrm
M}_T)\right)$.

Let $V_\ell$ (resp. $V_r$) be the space of left (resp. right)
continuous $\mathbb{R}$-valued functions with bounded variations. We put a
superscript $T$ to specify the functions on $[0,T]$. There is a
bijective correspondence between $V_r^T$ and ${\mathrm M}_T$ :

- for any $v \in V_r^T$, there exists a unique $\mu \in {\mathrm
M}_T$ such that $v = \mu([0, \cdot])$; we denote it by $\dot v$ ,

- for any $\mu \in {\mathrm M}_T$,  $v = \mu([0, \cdot])$ stands in
$V_r$.

For $v \in D$,  let $\dot v = \dot v_a + \dot v_s$ be the Lebesgue
decomposition of the measure $\dot v \in {\mathrm M}([0,1))$ in
absolutely continuous and singular parts with respect to the
Lebesgue measure and let $\mu$ be any bounded positive measure
dominating $\dot v_s$.

For $A \subset [0, 1)$ and $v \in D$  let \ben \label{intact0} I_A
(v) = \int_A L_a\Big(t, \frac{d\dot v_a}{dt} (t)\Big) dt + \int_A
L_s \Big(t, \frac{d\dot v_s}{d\mu} (t)\Big) d\mu(t)\ \ \ if \ \ v
\in V_r\,, \een and  $I_A (v) = \infty$ if $v \in D \setminus V_r$, where functions $L_a(t,x)$ and $L_s(t,x)$ will be defined later for each of the ensembles of interest. 

\subsubsection{Uniform Gram ensemble}
\label{subsubG} For the following statement, we need some notation.
Let ${\bf H}$ be the entropy function : \be{\bf H}(x|p)
=\displaystyle x \log \frac{x}{p} + (1-x) \log\frac{1-x}{1-p}\,,\ee
and put\footnote{we set $\delta(y |A) = 0$ if $y\in A$ and $= \infty$ if $y \notin A$}
 \ben\nonumber
 L_a^G (t, y) &=& {\bf H}(1-t | e^{y}) \ \delta(y |
(-\infty, 0))\,,\\\label{LaLs}  L_s^G (t , y) &=& -(1-t)y \ \delta(y
| (-\infty, 0)) \,. \een
\begin{theorem}
\label{LDPH} The sequence $\{ n^{-1}\!\
 \log \Delta_n^{G, \beta}(t) , \ t \in [0,1)  \}_n$ satisfies a LDP in $(D, \sigma(D, {\mathrm M}_<))$ in the scale $2\beta^{-1}n^{-2}$
with good rate function $I_{[0,1)}^G$.
\end{theorem}

That means, roughly speaking, that
$$\mathbb{P}( \log \Delta^{G, \beta}_n \simeq nv) \approx \e^{-\frac{\beta n^2}{2} I_{[0,1)}^G (v)}\,.$$
The proof, in Section \ref{7.4G}, needs several steps. Let
$\Theta^G_n = n^{-1}\!\ \log \Delta^{G, \beta}_n$, so that
\ben\label{rmesG}\dot \Theta^G_n = \frac{1}{n}\sum_{j=1}^n \Big(\log
\rho_{n,j}^{G, \beta}\Big)\!\ \delta_{j/n}\,.\een
 First we show that $\{ \dot \Theta^G_n\}$  satisfies a LDP in
${\mathrm M}_T$ equipped with the topology $\sigma({\mathrm M}_T,
V_\ell)$. Then we carry the LDP to $\left(D_T , \sigma(D_T, {\mathrm
M}_T)\right)$ with good rate function: \ben \label{42} I_{[0,T]}^G
(v) = \int_{[0,T]} L_a^G\Big(t, \frac{d\dot v_a}{dt} (t)\Big) dt +
\int_{[0,T]} L_s^G \Big(t, \frac{d\dot v_s}{d\mu} (t)\Big)
d\mu(t)\,. \een
 To end the proof we apply the
Dawson-G\"artner theorem on projective limits (\cite{DZ} Theorem
4.6.1, see also  \cite{leo1} Proposition A2).

 Let
us note that $I_{[0,T]}^G(v)$ vanishes only when $v$ satisfies
(essentially) \ben \label{vanish} \frac{d\dot v_a}{dt}(t) = \log
(1-t)\ \ , \ \ \frac{d\dot v_s}{d\mu}(t) = 0\,,\een
 i.e. for $v
(t) = - {\mathcal J}(1-t)$, which is consistent with the result
(\ref{theocvps}).
\medskip

The LDP for marginals is given in the following theorem, where a
rate function with affine part appears.
\begin{theorem}
\label{margH} For every $T < 1$, the sequence $\left\{n^{-1}\!\ \log
\Delta_n^{G, \beta}(T)\right\}_n$  satisfies a LDP in $\mathbb{R}$ in the
scale $2\beta^{-1}n^{-2}$ with good rate function denoted by \ben
\label{opt1G} I_T^{G} (\xi) = \inf \{I^G_{[0,T]}(v) \ ; \ v(T) =
\xi\}\,. \een
\begin{enumerate}
\item
 If $\xi \in  [-T, 0)$ the equation
\ben \label{=xi1H} {\mathcal J}(1+\theta) -{\mathcal J}(1-T+\theta)
- T \log (1+\theta) =\xi\,, \een has a unique solution, and we have
\ben \label{idetH}
I_T^{G}(\xi) = \theta\xi &+& T {\mathcal J} (1+\theta)\\
\nonumber &+& \left(F(1)-F(1-T) -F(1+\theta) +F(1-T+\theta)\right)
\,. \een
\item
 If $\xi < - T$, we have
\ben \label{mauvH} I_T^{G}(\xi) = I_T^{G}(-T) - (1-T)(\xi + T)\,.
\een
\item
If $\xi \geq 0$, $I_T^{G}(\xi) = \infty$.
\end{enumerate}
\end{theorem}

\subsubsection{Laguerre ensemble} Let \ben \label{identif} \nonumber L_a^L(t,y) &=&
(e^y -1)  -  (1-t)y +{\mathcal J}(1-t)\\
L_s^L (t,y) &=&  - (1-t) y \ \delta(y | (-\infty , 0)) \,. \een
\label{subsubL}
\begin{theorem}
\label{LDPwish} The sequence $\{n^{-1}\!\ \log \Delta^{L,
\beta}_{n}(t) , t \in [0, 1)\}_n$ satisfies a LDP in $(D,
\sigma(D, {\mathrm M}_<))$, in the scale $2\beta^{-1}n^{-2}$ with
good rate function $I^{L}_{[0,1)}$.
\end{theorem}
That means, roughly speaking, that
$$\mathbb{P}( \log \Delta^L_n \simeq nv) \approx \e^{-\frac{\beta n^2}{2} I_{[0,1)}^L (v)}\,.$$
The proof uses the above result for the Uniform Gram process and the
beta-gamma algebra.

Let us note that $I_{[0,T]}^L(v)$ vanishes only when $v$ satisfies
(\ref{vanish}) (again) i.e.  for $v (t) = - {\mathcal J}(1-t)$,
which is consistent with the result (\ref{cvpW}).
\medskip

The LDP for marginals is given in the following theorem.
\begin{theorem}
\label{margW} For every $T < 1$, the sequence $\{n^{-1}\!\ \log
\Delta^{L, \beta}_{n}(T) \}_n$
 satisfies a LDP
in $\mathbb{R}$ in the scale  $2\beta^{-1}n^{-2}$ with good rate function
denoted by $I_T^L$. \ben \label{opt1L} I_T^{L} (\xi) = \inf
\{I^L_{[0,T]}(v) \ ; \ v(T) = \xi\}\,. \een
\begin{enumerate}
\item If $\xi \geq \xi_T := {\mathcal J}(T) - 1$ the equation
\ben
\label{=xi1W}
 {\mathcal J}(1+\theta) -{\mathcal J}(1-T+\theta) =\xi\,.
\een
has a unique solution, and we have
\ben
\label{idetW}
I_T^L(\xi) =  \theta\xi
+ F(1)-F(1-T) -F(1+\theta) +F(1-T+\theta)
\,.
\een
\item
 If $\xi < \xi_T$
, we have
\ben
\label{mauvW}
I_T^L(\xi) = I_T^L(T) + (1-T) (\xi_T-\xi)\,.
\een
\end{enumerate}
\end{theorem}
\medskip

\subsubsection{Jacobi ensemble}
In this subsection, the endpoint $1$ of  subsections \ref{subsubG}
and \ref{subsubL} is  replaced by $\tau_1$.

 Let, for $t < \tau_1$, \ben\nonumber L_a^J(t,y)
&=& \left(\tau_1 + \tau_2
 -t\right){\bf H}\Big(\frac{\tau_1 -t
  }{\tau_1 + \tau_2
   -t}\!\ \Big|\  e^y\Big)
\\\label{LaJ}L_s^J(t,y) &=& -(\tau_1 -t)y \ \hbox{if} \ y < 0\,. \een
\begin{theorem}
\label{LDPpath}
The sequence $\{
 n^{-1}\log \Delta_n^{J, \beta}(t) , \ t \in [0, \tau_1)  \}_n$ satisfies a LDP
 in $(D, \sigma(D, {\mathrm M}_<))$ in the scale $2\beta^{-1}n^{-2}$
with good rate function $I_{[0, \tau_1)}$.
\end{theorem}
That means, roughly speaking, that
$$\mathbb{P}( \log \Delta^{J, \beta}_n \simeq nv) \approx \e^{-\frac{\beta n^2}{2} I_{[0,1)}^J (v)}\,.$$
Let us note that $I_{[0,T]}^J(v)$ vanishes only when $v$ satisfies
(essentially)
\[\frac{d\dot v_a}{dt}(t) = \log\frac{\tau_1 -t}{\tau_1 + \tau_2
-t}\ \ , \ \ \frac{d\dot v_s}{d\mu}(t) = 0\,, \]i.e. for $v (t) =
{\mathrm E}(\tau_1, \tau_2 ,t)$, which is consistent with the result
(\ref{esplim}).

The LDP for marginals is given in the following theorem.
\begin{theorem}
\label{margJ} Let $T \in [0,  \tau_1)$, and  $\xi_T^J = {\mathcal
J}(\tau_2) + {\mathcal J}(T) -{\mathcal J}(T + \tau_2) -1$.
\begin{enumerate}
\item The sequence $\{n^{-1}\!\  \log \Delta_n^{J, \beta}(T)\}_n$  satisfies a LDP in $\mathbb{R}$
in the scale  $2\beta^{-1}n^{-2}$ with good rate function $I_T^{J}$
where \ben \label{opt1J} I_{T}^{J} (\xi) := \inf \{I^J_{[0,T]}(v) \
; \ v(T) = \xi\}\,. \een
\item
 If $\xi \in [\xi_T^J , 0)$,
 the equation \ben \label{crucial0} {\mathcal E}(\theta+\tau_1,
\tau_2 , T)= \xi \een has a unique solution $\theta \geq T -\tau_1
$, and we have \ben \nonumber I_{T}^{J}(\xi) = \theta\xi &-&
\left[F(\theta+\tau_1)- F(\theta+\tau_1 -T) \right] -\left[ F(\tau_1 + \tau_2) - F(\tau_1 + \tau_2 -t)\right]\\
\nonumber
&+&
 \left[F(\tau_1)- F(\tau_1 -T)\right]\\
\label{idetJ}
&+& \left[F(\theta+\tau_1+\tau_2)- F(\theta+\tau_1 +\tau_2 -T)\right]
\,.
\een
\item
 If $\xi <\xi_T^J$,
  we have
\ben
\label{mauvJ}
I_T^{J}(\xi) = I_T^{J}(\xi_T^J) + (\xi_T^J - \xi)(\tau_1 - T)\,.
\een
\item
If $\xi \geq 0$, then $I_T^{J}(\xi) = \infty$.
\end{enumerate}
\end{theorem}
\section{Connections with the spectral method}
\label{contrac}
The logarithm of the determinant of a non singular matrix is a linear statistic of the empirical distribution of its eigenvalues, so that we may compare the above result with those obtained by this spectral approach.
\subsection{Laguerre/Wishart} \label{MP}
We start with \[\frac{1}{n}\log \Delta_{n,r}^{L, \beta}  =
\frac{r}{n}\int (\log x) \ d\mu_{n, r}(x)\]
where $\mu_{n,r}$ the so called empirical spectral distribution (ESD) is
\ben \label{defesdf} \mu_{n,r} = \frac{1}{r}
\sum_{k=1}^r \delta_{\lambda_k} \een

For $c > 0$ and $\sigma > 0$, let $\pi_{\sigma^2}^c$ be the
 distribution on $\mathbb {R}$ defined by \ben \label{defpi}
\pi_{\sigma^2}^c (dx)= (1-c^{-1})_{_+} \delta_0 (dx) +
\frac{\left((x - \sigma^2 a(c))(\sigma^2 b(c) - x)\right)_+^{1/2}}
{2\pi \sigma^2c x}\ dx\,, \een where $\delta_0$ is the Dirac mass in
$0$, $x_{_+} = \max (x,0)$ and \ben \label{aetb} a(c) = (1 -
\sqrt{c})^2 \ , \ \ b(c) = (1 + \sqrt{c})^2\,. \een It is called the
Mar\v cenko-Pastur distribution with ratio index $c$ and scale index
$\sigma^2$ (\cite{Baimethodo} p.621).

It is well known (\cite{MarchPastur1}, 
\cite{Baimethodo} section 2.1.2 for the cases $\beta=1$ and $\beta=2$) that as $n,r \rightarrow \infty$
with $r/n \rightarrow T \in (0, \infty)$, the family of ESD
$(\mu_{n,r})$
 converges a.s. weakly
to
 $\pi_1^T$. If we replace the
common law ${\mathcal N}(0 , 1)$ by ${\mathcal N}(0, \sigma^2)$ then
the limiting distribution is $\pi_{\sigma^2}^T$.

To conclude that
\ben \label{logW} \lim_n \int (\log x) \ d\mu_{n,
r} (x) = \int(\log x) \!\ d\pi_1^T(x)\,,\een an additional control
is necessary, since $x \mapsto \log x$ is not bounded.

Actually,
 the largest and the smallest eigenvalue
converge a.s. to $b(T) < \infty$ and  $a(T) > 0$, respectively. For
comments on these results and references, one may consult
\cite{Baimethodo} sections 2.1.2 and 2.2.2., (see also \cite{johnstone}).
In our context, this implies easily that a.s.
\ben \label{esp}
 \frac{1}{n}\log \Delta^{L, \beta}_{n,r}  = \frac{r}{n}\int (\log x) \ d\mu_{n,r}(x)
 \rightarrow T \int \log x  \ d\pi_1^T(x)
\een
 Moreover, it is known
(\cite{jonsson1} p.31 and 
\cite{BaiSilver2} p.596-597) that : \ben \label{often}\nonumber
T \int (\log x)  \ d\pi_1^T(x) &=& \int_{a(T)}^{b(T)}   \frac{\log x}{2\pi x}\ \sqrt{(x-a(T)) (b(T)-x)}\ dx\\
&=& (T-1) \log(1-T) -T = -{\mathcal J}(1-T) \een which implies that
claim (\ref{esp}) is  consistent with
(\ref{cvpW}).

Recently,
\cite{BaiSilver2} proved a CLT for linear statistics of sample covariance matrices (non necessarily Gaussian),
with the meaningful example of determinants. They consider the real and complex case, 
and their results (Theorem 1.1 ii) and iii) are  consistent with the marginal
version of (\ref{etantW}).
It is likely that $\beta=4$ can also be handled under their assumptions.
\medskip

Let us end with the large deviations.
\cite{hiai1}, (see also \cite{HiaiP} section 5.5)
proved\footnote{Their $\beta$ is our $\beta'$.} that if
$n\rightarrow\infty$ and $r/n \rightarrow T < 1$, then
$\{\mu_{n,r}\}$ satisfies a LDP in ${\mathrm M}_1([0, \infty))$ in
the scale $2\beta^{-1}n^{-2}$ with some explicit good rate function
$I_T^{spL}$
 given below in (\ref{HP1}, \ref{HP2},
\ref{HP3}). If the contraction $\mu \mapsto \int \log x \ d\mu(x)$
 were continuous, we would claim that
$\{n^{-1} \log \Delta^L_{n,\lfloor nT\rfloor}\}_n$ satisfies a LDP
in $\mathbb {R}$ in the same scale, with good rate function \ben \label{58}
\wt I^{L}_T (\xi) = \inf \left\{\wt I^{spL}_T (\mu) \ ; \ T\int \log
x \ d\mu(x) = \xi \right\}\,. \een Actually, \ben \label{HP1}
I^{spL}_T (\mu) = -T^2 \Sigma(\mu) + T \int\left(x - (1-T) \log
x\right) d\mu(x) + 2 B(T) \een where \ben \label{HP2} \Sigma(\mu) :=
\int\!\!\int \log|x-y|\ d\mu(x) d\mu(y) \een is the so-called
logarithmic entropy and for $T \in (0,1)$ \ben \label{HP3} 2 B(T) =
-\frac{1}{2}\left( 3T -T^2\log T + (1-T)^2 \log(1-T)\right)\,. \een

We do not know if the contraction $\mu \mapsto \int \log x \
d\mu(x)$ does work, although not continuous. However we will prove
the following result, where for $u \in \mathbb {R}$ we put \ben \label{calA}
{\mathcal A}(u) = \{\mu : \int (\log x) \ d\mu(x) =u\}\,. \een
\begin{prop}
\label{infIW} For $\xi\geq \xi_T$ and $\theta$ solution of
(\ref{=xi1W}), let  $\sigma^2 = 1 +\theta$. Then the infimum of $
I^{spL}_T (\mu)$ over ${\mathcal A}(\xi/T)$ is uniquely achieved for
$\pi_{\sigma^2}^{T/\sigma^2}$
 and
\ben \label{defxiL} I_T^L(\xi) =  I^{spL}_T
(\pi_{\sigma^2}^{T/\sigma^2})
 = \inf \{ I^{spL}_T (\mu) ; \ \mu \in {\mathcal A}(\xi/T)\}\,.
\een
\end{prop}

\begin{remark}
\begin{enumerate}
\item The endpoint is $\xi_T = {\mathcal J}(T) - 1$, with
  $\sigma^2 = T$.

\item For $\xi < \xi_T$ we do not know what happens.
We can imagine that the infimum  in (\ref{defxiL}) has a solution in some extended space.
\end{enumerate}
\end{remark}
\subsection{Uniform Gram}
Let $\widetilde\lambda_k, k=1, \dots , r$  be the (real)
eigenvalues of $G$ in the Uniform Gram ensemble, and set
\ben \label{defesdftilde} \widetilde\mu_{n,r}
= \frac{1}{r} \sum_{k=1}^r \delta_{\widetilde\lambda_k}\,.
\een
For $\beta =1$, \cite{Spince1} 
proved that, as $n\rightarrow \infty$ and $r/n \rightarrow T \in (0,
\infty)$, the family $(\widetilde\mu_{n,r})$ converges a.s. to
$\pi_1^T$.
More recently \cite{Jiang2} proved that 
 the same result holds true in a complex Gram ensemble not necessarily uniform.
Again, like in Section \ref{MP},
we may write
$$\frac{1}{n}\log \Delta_{n,r}^{G, \beta}  = \frac{r}{n}\int (\log x) \ d\widetilde\mu_{n, r}(x)$$
and use the weak convergence of $\widetilde\mu_{n, r}$ towards $\pi_1^T$. Recently,
\cite{Jiang2} proved  that the largest and the
smallest eigenvalue  converge a.s.  as $r/n \rightarrow T < 1$ to
$b(T)< \infty$ and $a(T)> 0$ respectively. So, we have \ben
\label{logG} \lim_n \int (\log x) \ d\widetilde\mu_{n, r} (x) = \int
(\log x) \!\ d\pi_1^T(x)\,. \een In view of (\ref{often}), this
matches with the result (\ref{theocvps}).

No result on fluctuations or large deviations seems to be known on $\widetilde\mu_{n,r}$.
\subsection{Jacobi} In the matrix models ($\beta =
1,2$ or $4$), take $r \leq n_1$ and let $\lambda_k , k=1, \dots , r$
be the  eigenvalues of ${\mathcal Z}_{n_1,n_2,r}$ (they are real
nonnegative). The ESD is
\be \nu_{n_1,n_2,r} = \frac{1}{r}
\sum_{k=1}^r \delta_{\lambda_k } \,. \ee When $n_2 \leq r \leq n_1$
we have
$$\nu_{n_1,n_2,r} = \frac{n_2}{r}\mu_{n_1,n_2,r} + \left(1-\frac{n_2}{r}\right)\delta_1\,,$$
where $\mu_{n_1,n_2,r}$ is the ESD built with eigenvalues different
from $1$. We can write in all cases
\ben
\label{logdet} \log \Delta_{n,r}^{J, \beta} =   \min(r, n_2)
\int (\log x)\ \! \mu_{n_1,n_2,r} (dx)\,.
\een
It is then possible to carry asymptotical results of this empirical distribution to $\log \Delta_{n,r}^{J,\beta}$.
 
\cite{capitaine} studied the complex case in
the asymptotical regime $n_1/r \rightarrow u', \ n_2/ r \rightarrow
v'$ with $u' + v' \geq 1$. They prove\footnote{They use the notation
$\alpha$ and $\beta$ but we change not to confuse with $\beta$
already defined.}
  that $\mathbb{E}\nu_{n_1,n_2,r}$ converges (in moments hence) in
  distribution. To give the expression of the limiting distribution, which we denote
$\hbox{CC}_{u', v'}$ and to compare with known results in some other
contexts with coherent notation, we will use in the following, four
functions :

\noindent for $(b,c) \in (0,1)\times(0,1)$ we put \ben
\label{defsig} \sigma_\pm (b,c) = \frac{1}{2}\left[1 + \sqrt{bc} \pm
\sqrt{(1- b)(1- c)}\right]\,, \een and for $(x,y) \in
(0,1)\times(0,1)$ \ben \nonumber
a_\pm(x,y) &=& (1 -x-y + 2xy) \pm 2 \sqrt{x(1-x)y(1-y)}\\
\label{defla} &=& \left(\sqrt{(1-x)(1-y)} \pm \sqrt{xy}\right)^2\,.
\een The mappings $\sigma_\pm$ and $a_\pm$ are inverse in the
following sense : \ben\label{lets}\{(b,c) : 0 < b < c < 1\}
\xrightleftharpoons[(a_-, a_+)]{(\sigma_- , \sigma_+)} \{(x,y) : 0<
x < y<1 \ \hbox{and} \ x+y > 1\}\een
 For $0 < a_- < a_+ < 1$, let $\pi_{a_-,a_+}$
 be the distribution on $\mathbb {R}$ defined by
\ben \label{defpij} \pi_{a_-,a_+}(dx)= C_{a_-,a_+} \frac{\sqrt{(x -
a_-)(a_+ - x)}}{2\pi x(1-x)}\!\ \mathbf{1}_{[a_-,a_+]}(x)\ dx\,, \een
where $C_{a_-,a_+}$ is the normalization constant. Since we found
some mistakes in the literature, let us compute explicitly the
constant $C_{a_-, a_+}$. From the obvious decomposition
$$\frac{1}{x(1-x)} = \frac{1}{x}+ \frac{1}{1-x}$$
we get 
\be (C_{a_-,a_+})^{-1} = I(a_-, a_+) + I(1-a_+, 1-a_-)\ee
where, for $0 < u <v$
\be
I(u,v) = \int_u^v \frac{\sqrt{(x -
u)(v - x)}}{2\pi x}\!\ dx 
\ee
This last integral could be calculated by elementary method, but it is shorter to connect it with the Mar\v cenko-Pastur distribution. Taking 
\[\sigma^2= \frac{\sqrt{v} + \sqrt{u}}{4}\ \ , \ \ \sqrt{c} = \frac{\sqrt{v} - \sqrt{u}}{\sqrt{v} + \sqrt{u}}\]
the simple fact that $\pi^c_{\sigma^2}$, given in (\ref{defpi}), is a probability distribution yields
\be I(u,v) = \frac{\left(\sqrt{v} -\sqrt{u}\right)^2}{4} \ \ \ (0<u<v)\,.\ee
Finally, we get:
\begin{equation}
\label{cste} (C_{a_-,a_+})^{-1}  =\frac{1}{2}\left[1 - \sqrt{a_-a_+}
- \sqrt{(1-a_-)(1-a_+)}\right]\,.
\end{equation}
 The distribution $\hbox{CC}_{u', v'}$ is then (recall $u' + v'
\geq 1$) : \begin{eqnarray} \label{CCdef}
 \hbox{CC}_{u', v'} := (1- u')^+\delta_0 &+&
(1-v')^+ \delta_1 \\ \nonumber&+&\left[1 - (1- u')^+ - (1-v')^+
\right]\pi_{a_-,a_+}\,, \end{eqnarray} where \begin{equation} \label{lambdalpha} (a_-,a_+)
= a_\pm\left(\frac{u'}{u' + v'} , 1 - \frac{1}{u' + v'}\right)\,.
\end{equation}

\begin{remark}
The case $(v' < 1)$ corresponds to $r > n_2$, the second matrix
$W_2$ is singular  and the case  $(v' \geq 1)$ corresponds to $r
\leq n_2$, the second matrix is non-singular.
\end{remark}

For particular values of the parameters and up to an affine change
to make the distribution symmetric, the distribution $\pi_{a_-,a_+}$
was introduced by \cite{Kesten} as limit distribution for random walks on
some classical groups. It was (independently)
introduced by  \cite{McKay}  as a limit distribution in a graph
problem. It is sometimes called the generalized McKay distribution.
Some important  connections are in Section \ref{appendixMK}.

For the LLN, the same remarks as above are relevant. Let us recall the notation
\[r \leq n_1\ , \ n\rightarrow \infty \ , \ \frac{r}{n}\rightarrow T \ , \  \frac{n_1}{n}\rightarrow \tau_1 \ , \ \frac{n_2}{n}\rightarrow \tau_2\ , \ u'= \frac{\tau_1}{T}\ , \ v' = \frac{\tau_2}{T}\,.\]
The weak
convergence of the ESD
 (\cite{capitaine}) and the control on the extremal eigenvalues (\cite{Led04},
  \cite{Collins} and references therein), yield, if $u' \geq 1$
\begin{eqnarray}\nonumber \lim_n \frac{1}{r} \log \Delta_n^{J, \beta}(T)
 &=& \int (\log
x) \!\ \hbox{CC}_{u', v'}(dx)\\
\label{interm}
&=& \hbox{min}(v',1) \int (\log x) \!\
 \pi_{a_-, a_+}(dx) \end{eqnarray}
where $a_\pm$ are in (\ref{lambdalpha}).
 Nevertheless a computation of this
integral by elementary methods is not so easy. After some attempts,
we choose to consider the above result as an indirect way to compute
this integral and we obtain the following result.
\begin{prop}
\label{intlogJ}For $0 < a_- < a_+ < 1$,  \ben \label{magic}
\int (\log x) \ \pi_{a_-,a_+}(dx) =\ \ \ \ \ \ \ \ \ \ \ \ \ \ \ \ \
\ \ \ \ \ \ \ \ \ \ \ \ \ \ \ \ \ \ \ \ \ \ \ \ \ \ \ \ \ \ \ \ \ \
\ \ \ \ \ \  \\\nonumber \ \ \ \ \ \ = \frac{\sigma_+ \log \sigma_+
+ \sigma_- \log \sigma_- -(\sigma_+ + \sigma_- -1)\log(\sigma_+ +
\sigma_- -1) }{1- \sigma_+} \een where $\sigma_\pm$ are specified by
(\ref{defsig}).
\end{prop}
\prf
From (\ref{esplim}), \be \lim_n
\frac{1}{r} \log \Delta_n^{J, \beta}(T)
 = \frac{1}{T} \lim_n
\frac{1}{n} \log  \Delta_n^{J, \beta}(T)
= \frac{1}{T}\!\
{\mathcal E}(\tau_1 , \tau_2 ,T) = {\mathcal E}(u' , v' ,1)\,,\ee
where for the last equality
we noticed that
${\mathcal E}$ is homogenous.
 With the help of (\ref{interm}) we get  \ben
\label{E=int} \hbox{min}(v',1) \int (\log x) \pi_{a_-, a_+}(dx) =
{\mathcal E}(u' , v' ,1)\,. \een From (\ref{lets}) we see that if
$u' \geq 1$ then
\[
\{\sigma_- , \sigma_+\} = \left\{\frac{u'}{u' + v'} ,
\frac{u'+v'-1}{u' + v'}\right\}
\]
We have two cases. When $v' > 1$,
$$\sigma_- = \frac{u'}{u' + v'} \ \ , \ \ \sigma_+ = \frac{u'+v'-1}{u' + v'}\,,$$
 so that
(\ref{E=int}) yields \ben \label{beta1} \int (\log x) \pi_{a_-,
a_+}(dx) = {\mathcal E}(u' , v' ,1) = {\mathcal
E}\Big(\frac{\sigma_-}{1-\sigma_+}, \frac{1- \sigma_-}{1-\sigma_+} ,
1\Big) \een
 When  $v' < 1$,  $$\sigma_+ = \frac{u'}{u'
+ v'}\ \ , \ \ \sigma_- = \frac{u'+v'-1}{u' + v'}\,,$$ so that
(\ref{E=int}) yields \ben \label{beta2} \int (\log x) \pi_{a_-,
a_+}(dx) = \frac{1}{v'}{\mathcal E}(u' , v' ,1) = \frac{1-
\sigma_-}{1-\sigma_+}\ {\mathcal E}\Big(\frac{\sigma_+}{1-\sigma_-},
\frac{1- \sigma_+}{1-\sigma_-} , 1\Big) \een and together
(\ref{beta1}-\ref{beta2}) provide (\ref{magic}). This ends the
proof.

Let us remark that the first case above $(v' < 1)$ corresponds to $T
> \tau_2$ (i.e. $r > n_2$, the second matrix $W_2$ is singular)  and
the second one $(v' \geq 1)$ corresponds to $T \leq \tau_2$ (i.e.
$r \leq n_2$, the second matrix is non-singular). $\Box$
\bigskip

Let us end with the large deviations. In the complex case ($\beta=
2$),  \cite{hiai2} proved that if $n \rightarrow \infty,
n_1/n \rightarrow \tau_1, n_2/n \rightarrow \tau_2 > \tau_1, r/n
\rightarrow T < \tau_1$, then $\{\mu_{n_1,n_2,r}\}_n$ satisfies a
LDP in   ${\mathcal M}_1([0,1])$ the set of probability measures on
$[0,1]$ endowed with the weak convergence topology, in the scale
$n^{-2}$, with the good rate function \ben  \nonumber I^{spJ}_{T}
(\mu) &:=& - T^2\Sigma(\mu)- T \int_0^1\left((\tau_1 -T) \log x +
(\tau_2 -T) \log (1-x)\right) \!\ d\mu(x)\\\label{hpldp}
 &+&  T^2  B\Big(\frac{\tau_1 -T}{T}, \frac{\tau_2 -T}{T}\Big)\,,
\een where $B$ is defined in (\ref{b1}) (it is the limiting free energy).

A computation similar to p.10 of \cite{hiai2} gives the same result
for general $\beta$.
\begin{prop}
\label{taux} If $T < \tau_1 \leq \tau_2$, the family $\{\mu_{\lfloor
n\tau_1\rfloor,\lfloor n\tau_2\rfloor, \lfloor nT\rfloor}\}$
satisfies a LDP in  ${\mathcal M}_1([0,1])$ in the scale
$2\beta^{-1}n^{-2}$ and good rate function $I^{spJ}_{T}$.
\end{prop}
If the contraction $\mu \mapsto \int \log x \ d\mu(x)$ from the set
${\mathrm M}_1([0,1])$
 to $\mathbb{R}$
 were continuous, we would claim that
$\{n^{-1} \log\Delta_n^{J, \beta}(T)\}_n$ satisfies a LDP in $\mathbb {R}$
with good rate function $\widetilde I^{J}_T$ where \ben
\label{infim} \widetilde I^{J}_T (\xi) = \inf \left\{I^{spJ}_T (\mu)
\ ; \ \mu \in {\mathcal A}(\xi T^{-1}) \right\} \een with ${\mathcal
A}(u)$ as defined in (\ref{calA}).

Like in the Laguerre case we will prove the following result.
\begin{prop}
\label{infI} Let $T < \min (\tau_1, \tau_2)$, $\xi \in [\xi_T^J ,
0)$ and $\theta$ solution of (\ref{crucial0}). Then the infimum of
$I_T^{spJ} (\mu)$ over ${\mathcal A}(\xi T^{-1})$ is uniquely
achieved at $\mu = \pi_{\tilde a_-, \tilde a_+}$ where
\[
(\tilde a_-, \tilde a_+) = a_\pm(\tilde s_- , \tilde s_+)
\]
with
\ben
\label{stilde}
\tilde s_- = \frac{\tau_1 + \theta}{\tau_1 + \tau_2 +\theta}\ , \ \tilde s_+ = \frac{\tau_1 + \theta + \tau_2  -t}{\tau_1 + \tau_2 +\theta}\,,
\een
 and
\ben \label{defxiJ} I_T^J(\xi) = I_T^{spJ} (\pi_{\tilde a_-, \tilde
a_+})
 = \inf \{ I_T^{spJ} (\mu) ; \ \mu \in {\mathcal A}(\xi T^{-1})\}\,.
\een
\end{prop}
\begin{remark}
The endpoint is $\xi_T^J$, which corresponds
 to $\theta = T- \tau_1$,
i.e. $$\tilde a_- = 0 \ \ , \ \  \tilde a_+ = \frac{4\tau_2 T}{(\tau_2 + T)^2}\,.$$
For $\xi < \xi_T^J$ we do not know what happens.
We can imagine that the infimum  in (\ref{infim}) has a solution in some extended space.
\end{remark}
\begin{remark}
In the range $\tau_2 \leq T < \tau_1$ we have a similar result, exchanging $\tilde s_-$ and $\tilde s_+$ in (\ref{stilde}). We omit the details.
\end{remark}
\subsection{Extensions}
We already mentioned that in the Wishart and Gram models, limiting
results exist for marginals when we leave the Gaussian/Uniform world, in particular for  fluctuations in \cite{BaiSilver2}.

The Bartlett decomposition is not possible in the general case.
Nevertheless, a product formula for the determinant is well known
(see for example  Lemma 3.1 p.9  and formula 4.3 p.15 in  \cite{Zeit}), but nothing can be said about the distribution
of the components of the product in general.

Nevertheless, if the columns (or the rows) of the matrix $B$ are
i.i.d. and isotropic, the previous results extend easily.

Let us begin with the "column" case. The beta-gamma algebra allowed
us to pass from the Uniform Gram ensemble to the Wishart ensemble.
 The polar decomposition allows to obtain similar results as
for the Wishart ensemble under convenient assumptions on the radial
distribution.
 Let $\varepsilon_n = \log \Vert b_1 \Vert^2 - \log \mathbb{E} \Vert b_1 \Vert^2$ (remember that we omit the dimension index $n$).
To get convergence and fluctuations it is enough to assume
 \ben
 \label{hypr}
 n \mathbb{E} \varepsilon_n \rightarrow a_1 \ , \ n\hbox{Var}\, \varepsilon_n \rightarrow a_2 \ , \ n \mathbb{E} (\varepsilon_n - \mathbb{E} \varepsilon_n)^4 \rightarrow 0\,.
 \een
To get large deviations, it would be sufficient to assume that, for some
convenient functions $\varphi$, the
quantity $n^{-2} \sum_{k=1}^n \ \log \mathbb{E} \exp
\left(n\varphi(k/n)\varepsilon_n\right) $ has a limit .

 \cite{akhavi1} uses the uniform distribution in the unit
{\it ball}, so that the distribution of $\Vert b_1\Vert^2$ is
Beta$\left(n/2, 1\right)$ and (\ref{hypr}) is satisfied with $a_1 =
-2 , \ a_2 = 0$. The contribution of the radial part is then roughly
"deterministic" since $\mathbb{E} \Vert b_1\Vert^2$ is bounded.

In the "row" case, we can use the results of the "column" case since
the  eigenvalues of $BB'$ are (except $0$ with multiplicity $n-r$)
the same as those of $B'B$.

\section{Proofs of Theorems of Section \ref{LLNF}}
\subsection{Proof of Theorem \ref{first2}}
We will use  Mellin transforms and their first two derivatives at
$\theta = 0$. From the decomposition (\ref{sumG}) we have  \ben
\label{melg}\log \mathbb{E} |\Delta^{G, \beta}_{n,r}|^{\beta' \theta} =
\sum_{k=1}^r \Lambda_{n,k}^{G,\beta} (\theta)  \een with  \ben
\label{LambdaG} \Lambda_{n,k}^{G, \beta} (\theta) := \log \mathbb{E}
\left[\rho_{n,k}^{G, \beta}\right]^{\beta'\theta}\een and from
(\ref{lawrhoG})\ben  \label{LambdaG+} \Lambda_{n,k}^{G, \beta}
(\theta) = \ell\big(\beta' (n-k+ 1+\theta)\big)\! - \!
\ell\big(\beta' (n-k+1)\big)\! + \! \ell(\beta' n)\!- \!
\ell\left(\beta' (n +\theta)\right) \een where we set
\[\ell(x) = \log \Gamma(x)\,.\]
\par\underline{Proof of 1) and 2)} Differentiating once, we get
\be \mathbb{E} \log \Delta^{G, \beta}_{n,r} =  \sum_{j=1}^r \left[\Psi
\left(\beta'(n-j+1)\right) - \Psi \left(\beta' n\right)\right]\,,
\ee and from Binet formula (\ref{bin3}), \ben \label{tripleG} \mathbb{E} \log \Delta^{G,
\beta}_{n,r} =  \log \frac{(n)_r}{n^r} +
\frac{1}{\beta}\left(H_{n-r} - H_n\right)+\frac{r}{\beta n}  - \delta_{n,r}^1
\,. \een in which

1)  $(p)_r = p(p-1)\cdots (p-r+1)$ is the falling factorial

2) $H_0 = 0$ and $H_p = 1 + \frac{1}{2} + \cdots + \frac{1}{p}$ are
the harmonic numbers

3) the delta term is  
\ben \label{delta1H} 
\ \delta^1_{n,r} =  \int_0^\infty sf(s) \sum_{k=1}^r [e^{-\beta' (n-k +1)s} - e^{-ns}]
\ ds 
\,.\een Using  Binet formula (\ref{bin1}) twice, we have for $r <
n$  \be\log \frac{(n)_r}{n^r} &=& -\Big(n -r + \frac{1}{2}\Big) \log
\left(1 - \frac{r}{n}\right) -r -
\int_0^\infty f(s) [e^{-s(n-r)} - e^{-sn}] ds\\
 &=& -n{\mathcal J}\left(1-\frac{r}{n}\right) -\frac{1}{2}
\log \left(1-\frac{r}{n}\right) - \int_0^\infty f(s) [e^{-s(n-r)} - e^{-sn}] ds\,. \ee
For $r=n$ the Stirling formula gives \[\log
\frac{(n)_n}{n^n}= -n + \frac{1}{2}\log(2\pi n) + o(1)\,,\]
 The
harmonic contribution in (\ref{tripleG}) is \be H_{n-r} - H_n = \log
\left(1 - \frac{r}{n}\right) + o(1) \ee as soon as $n-r \rightarrow
\infty$. For $r=n$, we have $H_0 - H_n = -\log n - \gamma + o(1)\,.$
Applying the dominated convergence theorem and (\ref{bina}), we see that the delta contribution satisfies:
 \[\sup_{r\leq n} \delta_{n,r}^1  = \delta_{n,n}^1 \rightarrow \int_0^\infty \frac{sf(s)}{e^{\beta' s} -1}\ ds\,,\]
 and $\lim_n \delta_{n, \lfloor nt\rfloor}^1 = 0$ for $t < 1$. 
%
Gathering all these estimates, and applying again the dominated convergence theorem, we get (for $n-r \rightarrow \infty$) \be \E \log \Delta^{G,
\beta}_{n,r} = -n {\mathcal J}\left(1-\frac{r}{n}\right) +
\frac{r}{\beta n}+\left(\frac{1}{\beta}-\frac{1}{2}\right) \log
\left(1 - \frac{r}{n}\right) + o(1)\,, \ee
 and
for $r=n$ \be \E \log \Delta^{G, \beta}_{n,n} = -n -
\left(\frac{1}{\beta} - \frac{1}{2}\right) \log n + K_\beta^1 +
o(1)\,. \ee Moreover, for the supremum, we have \be
\sup_{r\leq n}\left|\E \log \Delta^{G, \beta}_{n,r}  - \log \frac{(n)_r}{n^r}\right| &=& O(\log n)\\
\sup_{r\leq n}\left|\log \frac{(n)_r}{n^r} + n{\mathcal J}\left(1 -
\frac{r}{n}\right)\right|&=& O(\log n) \ee so that (\ref{enough}),
(\ref{esptG}) and (\ref{esp1G}) are proved.

3) Taking logarithms in (\ref{melg}) and differentiating twice, we
get \be \hbox{Var} \log \Delta^{G, \beta}_{n,r} =  \sum_{j=1}^r
\Psi' \left(\beta'(n-j+1)\right) - \Psi'(\beta' n)\ee and owing to
(\ref{restepsi}) \be \hbox{Var} \log \Delta^{G, \beta}_{n,r} =
\frac{1}{\beta'} (H_n - H_{n-r}) -\frac{r}{\beta' n} + \varepsilon
\ee where \be
 |\varepsilon| \leq \sum_{n-r+1}^n \frac{2}{\beta'^2 j^2}
\ee Moreover \be \hbox{Var} \log \Delta^{G, \beta}_{n,n} &=&
\frac{1}{\beta'} H_n -\frac{1}{\beta'}+ 
\delta_n^2\ee
where
\[\delta_n^2 = \int_0^\infty s\left(sf(s) +
\frac{1}{2}\right) \sum_{k=1}^n [e^{-\beta' (n-k +1)s} - e^{-\beta'ns}] \ ds\,.\]
Applying again the dominated convergence theorem and (\ref{proprf}), we get
\[\lim_n \delta_n^2 = \int_0^\infty \frac{s\left(sf(s) +
\frac{1}{2}\right)}{e^{\beta's}-1}\ ds\]
 Using (\ref{proprf}) and dominated
convergence we deduce easily (\ref{vartG}) and (\ref{var1G}).

To prove 4), let us note that since ${\mathcal J}$ is uniformly
continuous on $[0,1]$ we have
$$\lim_n \sup_{t\in [0,1]} \left|{\mathcal J}\left(1 - \frac{\lfloor nt\rfloor}{n}\right) - {\mathcal J}(1-t)\right| = 0\,,$$
so that, owing to (\ref{enough}), it is enough to prove that in
probability \[ \sup_{1\leq p\leq n}
 \left|\log \Delta_{p , n}^{G, \beta}
- \mathbb{E}\log \Delta_{p , n}^{G, \beta}\right| = o(n)\,.\] Actually this
convergence is a consequence of  Doob inequality and of the variance
estimate Var $n^{-1} \Delta^{G, \beta}_{n, n} = O(n^{-2}\log n)$
coming from (\ref{var1G}). $\Box$

\subsection{Proof of Theorem \ref{DonskerH}}
Let us first note that, thanks to the estimations of expectations
in (\ref{esptG}) and (\ref{esp1G}), we can reduce the problem to the
centered process and centered variable : \[\delta_n (t) := \log
\Delta^{G, \beta}_n (t) - \mathbb{E}\!\ \log \Delta^{G, \beta}_n (t) \ \ , \
\ \widehat\delta_n = \delta_n (1)/ \sqrt{(2/\beta)\log n}\,.\]

\noindent 1) We have $\delta_n (t) = \sum_{k=1}^{\lfloor nt\rfloor}
\eta_{n,k}$ where \ben \label{defeta} \eta_{n,k} := ( \log \rho^{G,
\beta}_{k,n}) - \mathbb{E}( \log \rho^{G, \beta}_{k,n}),\ \  k \leq n\een is
a row-wise independent arrow. To prove (\ref{etantG}) it is enough
to prove the convergence in distribution in $\mathbb {D}([0,T])$, for every
$T < 1$, of $\delta_n$ to a centered Gaussian process with
independent increments, and variance $\int_0^t \Big(\sigma^{G,
\beta}(s)\Big)^2\ ds$. To this purpose we apply a version of the
Lindeberg-L\'evy-Lyapunov criteria (see 
\citet*{DDC} Volume II Theorem 7.4.28 , or  \citet*[Chap.~3~c]{JacShi}). For $t < 1$, from (\ref{vartG}) it is
enough to prove that \ben \label{ddc4G} \lim_n \sum_{k=1}^{\lfloor
nt\rfloor} \mathbb{E}\!\ (\eta_{n,k}^4) = 0\,. \een
 We have from definitions (\ref{defeta}) and (\ref{LambdaG})
 \ben
\label{4=1+3G} \beta'^4 \mathbb{E} (\eta_{n,k}^4)
 =
(\Lambda_{n,k}^{G, \beta}) ^{(4)} (0) + 3 [(\Lambda_{n,k}^{G,
\beta})^{(2)} (0)]^2\,, \een
 On the one hand, from expression (\ref{LambdaG+})
 \[(\Lambda_{n,k}^{G, \beta})^{(4)} (0) = \beta'^4
[\Psi^{(3)}(\beta'(n-k+1))- \Psi^{(3)}(\beta'n) ]\]
 and Binet estimates (\ref{polygamma}),
(\ref{restepsi}) for $q=4$ yield \ben \label{channel4G}
 \left|\sum_{k=1}^p
(\Lambda_{n,k}^{G, \beta}) ^{(4)} (0) - 6\beta' \sum_{k=1}^p\left[
\frac{1}{(n-k + 1)^3} - \frac{1}{n^3}\right] \right|
 \leq 6
  \sum_{k=1}^p \frac{1}{(n-k+1)^4}\,,
\een
 which, for $0< t < 1$ and $p=\lfloor nt\rfloor$ yields $\lim_n
\sum_{k=1}^{\lfloor nt\rfloor}
\Big(\Lambda_{n,k}^{G,\beta}\Big)^{(4)} (0)  =0\,.$ On the other
hand, \ben \label{supG} \sum_{k=1}^p [(\Lambda_{n,k}^{G, \beta}) ''
(0)]^2 \leq \Big(\sup_{j\leq p}(\Lambda_{n,j}^{G, \beta}) ''
(0)\Big) \sum_{k=1}^p (\Lambda_{n,k}^{G, \beta}) '' (0)\,. \een We
already know, from (\ref{vartG}) that
\[ \beta'^{-2} \sum_{k=1}^{\lfloor nt\rfloor} (\Lambda_{k,n}^{G, \beta}) '' (0)
= \hbox{Var} \log \Delta_n^{G, \beta}(t) \rightarrow \int_0^t
\Big(\sigma^{G, \beta}(s)\Big)^2 \ ds\,.\] Now since
$(\Lambda_{n,k}^{G, \beta})'' (0) = \beta'^2 [\Psi'(\beta'(n-k+1))-
\Psi'(\beta'n) ]$
and since  $\Psi'$ is non-increasing (see (\ref{polygamma})) we
obtain \be
 \sup_{j\leq \lfloor nt\rfloor}(\Lambda_{j,n}^{G, \beta}) ''
(0)\leq \beta'^2 \Psi'\left(\beta'(n-\lfloor nt\rfloor+1)\right)\,,
\ee and from (\ref{restepsi}) (again), this term  tends to $0$. We
just checked (\ref{ddc4G}), which proves that the sequence of
processes $\{\delta_n(t), t \in [0, 1)\}_n$ converges to a Gaussian
centered process  ${\mathcal W}$ with independent increments and the
convenient variance. It is now straightforward to get equation
(\ref{sdeG}).
\medskip

\noindent 2) When $t =1$, most of the sums studied above explode
when $n$ tends to infinity and we need a renormalisation. In fact,
for every $n$, the process $\left(\delta_n (t), t\in [0,1]\right)$
has independent increments. The conditional distribution of
$\delta_n (1)$, knowing
  $\delta_n (t_1)=\varepsilon_1, \dots , \delta_n (t_r) = \varepsilon_r$
for $t_1 < \cdots < t_r$ is the same as $\varepsilon_r +
\sum_{[nt_r]+1}^n \eta_{k,n}$. Formulae (\ref{vartG}) and
(\ref{var1G}) yield \ben \label{vlogG} \sum_{[nt_r]+1}^n \mathbb{E}
(\eta_{k,n}^2) = (2/\beta) \log n + O(1)\,. \een Actually we can
apply the Lindeberg theorem (with the criterion of Lyapunov) to the
triangular array of random variables $\widehat \eta_{k,n} =
\eta_{k,n}/ {\sqrt{(2/\beta)\log n}}$ with with $k = [nt_r]+1, \dots
, n$. It is enough to prove \ben \label{lyaG} \lim _n \sum_{k=1}^n
\mathbb{E} (\widehat \eta_{k,n}^4) = 0\,. \een We start again with the
decomposition (\ref{4=1+3G}). From the above estimate
(\ref{channel4G}), the sum
  $\sum_{k=1}^n (\Lambda_{n,k}^{G, \beta})^{(4)} (0)$
 is bounded.
In (\ref{supG}), we have
\[\sum_{k=1}^n (\Lambda_{n,k}^{G, \beta}) '' (0)= \beta'^{-2} \hbox{Var}\ \log \Delta_n^{G, \beta}(1) \]
which is  equivalent to $ 2\log n$ (see (\ref{var1G})) and the
supremum in (\ref{supG}) with $p=n$ is bounded.
 This yields
\[
\sum_{k=1}^n \mathbb{E} (\widehat \eta_{k,n}^4) = \beta'^2 (\log n)^{-2}
\sum_{k=1}^n \mathbb{E} (\eta_{k,n}^4) = O((\log n)^{-1} )
\]
which proves
 (\ref{lyaG}).

Then $\sum_{[nt_r]+1}^n  \widehat\eta_{k,n}$ converges in
distribution to ${\mathcal N}(0,1)$, and the same is true for the
conditional distribution of $\widehat\delta_n$ knowing $\delta_n
(t_1)=\varepsilon_1, \dots , \delta_n (t_r) = \varepsilon_r$. Since
the limiting distribution does not depend on $\varepsilon_1, \dots ,
\varepsilon_r$, we have proved that $\widehat\delta_n$ converges in
distribution to a random variable which is ${\mathcal N} (0,1)$ and
independent of ${\mathcal W}$. $\Box$
\subsection{Proof of Theorems \ref{ODD} and \ref{DonskerW}}
It is of course possible to follow the same schemes of proof.
Actually we prefer, at least for the beginning, exploit the
beta-gamma algebra and the fundamental relation (\ref{relproc}). So,
for instance \be \mathbb{E}\left[\varepsilon_k^{(n)}\right]^{\beta' \theta}
= \Big(\frac{1}{\beta' n}\Big)^{\beta' \theta}\frac{\Gamma(\beta'
(n+\theta))}{\Gamma(\beta'n)}\ee hence \ben \label{newcgf} \log \mathbb{E}
\left[\varepsilon_k^{(n)}\right]^{\beta' \theta} =
\ell\left(\beta'(\theta +n)\right) - \ell\left(\beta'n\right) -
\beta'\theta \log\left(\beta'n\right)\,, \een which provides
estimates for the expectation and the variance. Differentiating once
and taking $\theta=0$, we see that \be
 \mathbb{E} \log  \varepsilon_k^{(n)} =
 \Psi\left(\beta'n\right) - \log\left(\beta'n\right)
= -\frac{1}{\beta n}-\int_0^\infty e^{-s\beta'n} s f(s)\ ds \\=
-\frac{1}{\beta n} + O\Big(\frac{1}{n^{2}}\Big) \ee (see
(\ref{bin3}), (\ref{proprf})), which gives \ben \label{espaux}
\sup_{p \leq n}\left|\mathbb{E} S_{n,p} + \frac{p}{n\beta}\right| =
O\Big(\frac{1}{n}\Big)\,.\een Besides, differentiating
(\ref{newcgf}) twice and taking $\theta=0$ again, we have \be
\hbox{Var} \Big(\log \varepsilon_k^{(n)}\Big) =
\Psi'\left(\beta'n\right) = \frac{1}{\beta'n} +
O\Big(\frac{1}{n^{2}}\Big) \ee (see (\ref{restepsi})), which yields
\ben \label{varaux} \sup_{p \leq n} \left|\hbox{Var} \!\ S_{n,p} -
\frac{2p}{\beta n}\right| = O\Big(\frac{1}{n}\Big)\,. \een From
(\ref{espaux}) and $(\ref{varaux})$ it is easy to check (via a
fourth moment estimate) that
 $S_n$
 converges in distribution in $\mathbb {D}([0,1])$ to
 \[\Big(-(t/\beta) + \sqrt{2/\beta}\!\ \widetilde{\bf B}_t, \ t \in [0,1]\Big)\] where $\widetilde{\bf B}$
 is a Brownian motion independent of
 $(\Delta^{G,\beta}_n , n \in \mathbb {N})$. Finally
 the family of processes $\Delta_n^{L, \beta} = \Delta^{G,\beta}_n + S_n$ converges in distribution towards
 \[\Big( X^{G,\beta}_t -(t/\beta) + \sqrt{2/\beta}\!\ \widetilde{\bf B}_t , \ t \in [0,1)\Big)\,.\]
 It is a Gaussian process, whose drift and variance coefficients are
 \be  \d^{G, \beta} (t) -\frac{1}{\beta} = \left(\frac{1}{2}-\frac{1}{\beta}\right) \frac{1}{1-t} = \d^{L, \beta}(t)\
 \,, \ \
 \Big(\sigma^{G, \beta} (t)\Big)^2 + \frac{2}{\beta} = \Big(\sigma^{L, \beta}(t)\Big)^2\,.
 \ee
 which identify the process $X^{L, \beta}$.

 Besides, we have
 $$\widehat\eta^{L, \beta}_n (1) =  \widehat\eta^{G, \beta}_n (1) + \frac{S_n(1)}{\sqrt{2\log n}}\,,$$
 so that the convergence of $\widehat\eta^{L, \beta}_n (1)$ is clear. Moreover the independence properties seen
 in Theorem \ref{DonskerH} remain true. $\Box$

\subsection{Proof of Theorem \ref{thlln}}
Again, we could follow the same schemes as in the Gram section.
Actually we take again the benefit of beta-gamma algebra. Let us
delete the superscript $\beta$ for the sake of simplicity.

Let us recall the equality in law (\ref{ThetamoinsTheta})
\[\log \Delta_{n_1, r}^{L} \el \log \Delta_{n, r}^{J} +
\log \Delta_{n_1 + n_2, r}^{L}- r \log\frac{n_1}{n_1 + n_2}\]
with independence in the left hand side.

We deduce easily
\[
\mathbb{E} \log \Delta_{n, r}^{J} = \mathbb{E} \log \Delta_{n_1, r}^{L} - \mathbb{E} \log
\Delta_{n_1 + n_2, r}^{L} + r \log \frac{n_1}{n_1 + n_2}
\]
 and
 \[ \hbox{Var}\!\ \log \Delta_{n, r}^{J} = \hbox{Var}\!\ \log \Delta_{n_1, r}^{L} -  \hbox{Var}\!\ \log \Delta_{n_1 + n_2, r}^{L}\]
The results  are now straightforward. We let the proof to the
reader. We just note that since $r/n_1 \rightarrow t/\tau_1$ and
$r/ (n_1+ n_2) \rightarrow t/(\tau_1 + \tau_2)$ then \be\mathbb{E} \log
\Delta_{n_1, r}^{L}+ n_1 {\mathcal J} \left(1- \frac{r}{n_1}\right)
&\rightarrow& \int_0^{t/\tau_1} \d^{L, \beta}(s) ds \\
\mathbb{E} \log \Delta_{n_1 + n_2, r}^{L}  + (n_1 + n_2) {\mathcal J}
\left(1- \frac{r}{n_1+n_2}\right)&\rightarrow&\int_0^{t/(\tau_1+
\tau_2)} \d^{L, \beta}(s) ds\ee
 hence \[\d^{J,
\beta}(t) = \frac{1}{\tau_1}\d^{L}\left(\frac{t}{\tau_1}\right) -
\frac{1}{\tau_1 + \tau_2}\d^{L}\left(\frac{t}{\tau_1+
\tau_2}\right)\,.\] In the same vein \be \Big(\sigma^{J}(t)\Big)^2 =
\frac{1}{\tau_1}\Big(\sigma^{L}\left(\frac{t}{\tau_1}\right)\Big)^2
- \frac{1}{\tau_1 + \tau_2} \Big(\sigma^{L}\left(\frac{t}{\tau_1 +
\tau_2}\right)\Big)^2\,.\ee

\subsection{Proof of Theorem \ref{DonskerJ}}
\label{55}
 Again, it is possible to follow the classical scheme.
Instead let us look at the situation we are faced to. Put \ben
 \nonumber
&\ & U_n =\log
\Delta_{n_1,r}^{L} -\mathbb E \Delta_{n_1,r}^{L} \ , \
 V_n =  \log \Delta_{n_1 + n_2,r}^{L}  -\mathbb E \log \Delta_{n_1 + n_2,r}^{L}\,,\\
\label{u+v}
&\ & W_n 
=  \log \Delta_{n,r}^{J} -\mathbb E \log \Delta_{n,r}^{J}\,. \een
so that $U_n  = V_n + W_n $
with $U_n  \Rightarrow U$ and $V_n \Rightarrow V$, where $U$ and $V$ are Gaussian processes with
independent increments, and $V_n$ and $W_n$ are independent. Looking
for instance at characteristic functions, it is clear that $W_n$
converges in the sense of finite distributions to a Gaussian process
with independent increments. Its drift and variance are the
difference of the corresponding ones. Moreover, since $\{U_n\}_n$ and
$\{V_n \}_n$ are tight, $\{U_n- V_n\}_n$ is tight.

\section{Proofs of Theorems of Section \ref{LD}}
\subsection{Proof of Theorem \ref{LDPH}}
\label{7.4G} Recall the notation $\Theta^G_n = n^{-1} \log
\Delta_n^{G, \beta}$. As mentioned after the statement of the
theorem, we are going to prove at first the LDP for the restriction of
$ \dot\Theta^G_n$ to $[0,T]$, viewed as an element of ${\mathrm
M}_T$, in the scale $\beta'^{-1}n^{-2}$ with rate function \ben
\label{tauxHTG} \widetilde I_{[0,T]}^{G} (m) := \int_0^T
L_a^G\Big(t, \frac{dm_a}{dt} (t)\Big) dt + \int_0^T L_s^G \Big(t,
\frac{dm_s}{d\mu} (t)\Big) d\mu(t)\,. \een Let $V_\ell$ be the set
of functions from $[0,T]$ to $\mathbb{R}$ which are left continuous and
have bounded variation, and let $V_\ell^{*}$ be its topological dual
when $V_\ell$ is equipped with the uniform convergence topology.

 Actually  $ \dot\Theta^G_n  \in {\mathrm M}_T$ may be identified with an element of $V_\ell^{*}$
 (see \cite{leo1} Appendix B): owing to (\ref{rmesG}) its action on $\varphi \in V_\ell$ is given by
\be < \dot \Theta^G_n, \varphi> := \frac{1}{n} \sum_{k=1}^{\lfloor
nT\rfloor} \varphi(k/n)  \log \rho^{G, \beta}_{n,k}\,. \ee

The proof of Theorem \ref{LDPH} is based on the ideas of Baldi
theorem (\cite{DZ} p.157). The main tool is the normalized cumulant
generated function (n.c.g.f.) which here takes the form \ben
\label{defcgfnG} {\mathcal L}_{n, \lfloor nT\rfloor}^{G, \beta}
(\varphi) := \frac{1}{\beta' n^2}\log \mathbb{E}\left[ \exp \Big(\beta' n^2
<\dot \Theta^G_n , \varphi>\Big)\right]\,. \een Owing to
(\ref{cvcgfG}) we have \ben \label{cgfnG} {\mathcal L}_{n, \lfloor
nT\rfloor}^{G, \beta} (\varphi) = \frac{1}{\beta' n^2}
\sum_{k=1}^{\lfloor nT\rfloor} \Lambda_{n,k}^{G,
\beta}\left(n\varphi(k/n)\right) \een and from (\ref{LambdaG}) it is
finite iff $\varphi(k/n) > - (n-k+1)/n $ for every $1 \leq k \leq
\lfloor nT\rfloor$.

In Subsection \ref{cvncgf}, we  prove the convergence of this
sequence of n.c.g.f. for a large class of functions $\varphi$. It
will be sufficient, jointly to the variational formula given in
Subsection \ref{variaf} to get the upper bound for compact sets.
Then Subsection \ref{expotigh} is devoted to exponential tightness,
which allows to get the upper bound for closed sets. However, since
the limiting n.c.g.f. is not defined everywhere, the lower bound
(for open sets) is more delicate than in Baldi theorem. Actually a
careful study of exposed points as in  \cite{grz} is
managed in Subsection \ref{exposed}. We end the proof in Subsection
\ref{end}.

\subsubsection{Convergence of the n.c.g.f.}
\label{cvncgf} Let, for $t \in [0,1]$ and $\theta > - (1-t)$ \ben
\label{defg} g^G(t, \theta) := {\mathcal J}(1-t +\theta) -{\mathcal
J}(1-t)- {\mathcal J}(1 + \theta)\,. \een
\begin{lem}
\label{cvcgfG} If $\varphi \in V_\ell$ satisfies $\varphi(t) +1-t >
0$ for every $t\in (0,T]$, then \ben \label{98G} \lim_n {\mathcal
L}_{n, \lfloor nT\rfloor}^{G, \beta} (\varphi) =\Lambda_{[0,T]}^{G}
(\varphi) := \int_0^T g^G(t, \varphi(t)) \ dt \,. \een
\end{lem}

\prf The key point is a convergence of Riemann sums. From
(\ref{LambdaG}) and (\ref{bin1}) we have, for every $\theta > -
\frac{n-k+1}{n}$, \be \Lambda_{n,k}^{G, \beta} (n\theta) = \beta'
(n-k+n\theta) \log \Big( 1 - \frac{k}{n} + \theta +\frac{1}{n}\Big)
-
\beta'(n-k) \log \Big( 1 - \frac{k}{n} +\frac{1}{n}\Big) \\
-\beta'(n-1 + n\theta)\log (1+\theta) + R_{n,k} (\theta) \ee where
the quantity\be R_{n,k} (\theta) = \int_0^\infty f(s)
e^{-\beta's}\left[ \e^{-\beta'(n-k +n\theta)s} - \e^{-\beta'(n-k)s}
\right]\!\ ds\\ - \int_0^\infty f(s)\left[ \e^{-\beta'(n-1
+n\theta)s} -
 e^{-\beta's} \e^{-\beta'(n-1)s}\right]\!\
ds\ee is bounded by  $2 \int_0^\infty e^{-\beta's} f(s)\!\ ds$. If
we set \be
\Phi_n (t) &:=& (1-t + \varphi(t))\log \Big( 1 -t + \varphi(t) + \frac{1}{n}\Big)\\
\nonumber && - (1 -t) \log( 1-t+ \frac{1}{n})-\Big(1-\frac{1}{n} +
\varphi(t)\Big)\log (1+\varphi(t)) \ee then, making $\theta =
\varphi(k/n)$, and adding in $k$, we get from (\ref{cgfnG}) 
\be
\frac{1}{\beta'n^2}\Big({\mathcal L}_{n, \lfloor nt\rfloor}^G
(\varphi)- \sum_2^{\lfloor nt\rfloor} R_{n,k}(\varphi(k/n))\Big) =
\frac{1}{n}\sum_1^{\lfloor
nt\rfloor}\Phi_n\left(\frac{k}{n}\right)=\\ = \int_{1/n}^{\lfloor
nt\rfloor/n} \Phi_n \left(\frac{\lfloor ns\rfloor}{n}\right) ds +
\frac{1}{n} \Phi_n \left(\frac{\lfloor nt\rfloor}{n}\right)\,. \ee
On the one hand, since $\varphi$ is left continuous, $\lim_n
\Phi_n\left(\frac{\lfloor nt\rfloor}{n}\right) =   g(t, \varphi(t))$
for every $t\in [0,T]$. On the other hand the following double
inequality  holds true: \be \nonumber \Phi_n (t) &\geq& \left(1-t +
\varphi(t)\right) \log \left(1-t + \varphi(t)\right) - (1-t) \log
(2-t)
 \\&-& \left(1+ \varphi(t)\right)\log \left(1+ \varphi(t)\right) -|\log \left(1-t + \varphi(t)\right)|
\\
\nonumber
   \Phi_n (t)&\leq&
 \left(1-t + \varphi(t)\right) \log \left(2-t + \varphi(t)\right) - (1-t) \log (1-t) \\
 &-& \left( 1 + \varphi(t) \right) \log \left(1+ \varphi(t)\right) +|\log \left(1-t + \varphi(t)\right)|\,,
\ee and with our assumptions on $\varphi$,  these bounds are both
integrable. This allows to apply the dominated convergence theorem
which ends the proof of Lemma
 \ref{cvcgfH}. $\Box$
\bigskip
 
If there exists  $s < T$ such that  $\varphi(s) < - (1-s)$ then for
$n$ large enough, ${\mathcal L}_{n,\lfloor nT\rfloor}(\varphi) =
+\infty$ and we set $\Lambda_{[0,T]}^G(\varphi) = \infty$. In the
other cases we do not know what happens, but as in \cite{grz}, we
will study the exposed points. Before, we need another expression of
the dual of $\Lambda_{[0,T]}^G$.

\subsubsection{Variational formula}
\label{variaf} Let us define $\Lambda_{[0,T]}^G (\varphi) = +\infty$
if $\varphi$ does not satisfy the assumption of Lemma \ref{cvcgfG}.
The dual of $\Lambda_{[0,T]}^G$ is then \ben \label{dualeG1}
\Big(\Lambda_{[0,T]}^G\Big)^\star (\nu) = \sup_{\varphi \in V_\ell}
\left\{<\nu, \varphi> - \Lambda_{[0,T]}^G (\varphi) \right\} \een
for $\nu \in V_\ell^*$. Mimicking the method of 
\cite{leo1} p. 112-113, we get \ben \label{dualeG2}
\Big(\Lambda_{[0,T]}^G\Big)^\star (\nu) = \sup_{\varphi \in
{\mathcal C}} \left\{<\nu, \varphi> - \Lambda_{[0,T]}^G (\varphi)
\right\} \een where ${\mathcal C}$ is the set of continuous
functions from $[0,T]$ into $\mathbb {R}$ vanishing at $0$. Then we apply
Theorem 5 of  \cite{Rocky1} and get \be
\Big(\Lambda_{[0,T]}^G\Big)^\star (\nu) = \int_0^T g^\star\Big(t,
\frac{d\nu_a}{dt}\Big) \!\ dt + \int_0^T r\Big(t,
\frac{d\nu_s}{d\mu}(t)\Big)\!\ d\mu(t) \ee where \ben
\label{identifH} g^\star(t,y) &=& \sup_\lambda \left\{\lambda y -
g^G(t,\lambda)\delta(\lambda | (-1, \infty)) \right\}\,, \een and
$r$ is the
 recession function :
 $$r(t,y) = \lim_{\kappa \rightarrow \infty} \frac{g^\star (t,\kappa y)}{\kappa}\,.$$
 Actually, if $y < 0$, the supremum is achieved for
\ben \label{lambdaH} \lambda^G(t,y) := -\Big(1- \frac{t}{1-
e^{y}}\Big)\, \een and we have
 \ben
 \nonumber
g^\star(t,y)&=& \lambda^G(t,y) y -g^G \left(t, \lambda^G(t,y)\right)\\
\nonumber
&=& -y(1-t) +(1-t) \log(1-t) + t \log t -t \log (1-e^{y})\\
\label{gstar} &=&  {\bf H}\left(1-t | e^{y}\right) \,. \een If $y
\geq 0$, $g^\star(t,y) = \infty$.  The recession is now
 $ r(t,y )= - (1-t) y$ if $y \leq 0$, and $=\infty$ si $y >0$.
As a result \ben \label{gLa} g^\star (t,y) = L_a^G (t,y)\ \ , \ \
r(t,y) = L_s^G (t,y)\,. \een
 So we proved the identification $\Big(\Lambda_{[0,T]}^G\Big)^\star = \widetilde I^G_{[0,T]}$
(recall (\ref{tauxHTG})).

\subsubsection{Exponential tightness}
\label{expotigh} In this paragraph and in Section 6.2 we use the
function defined for $\theta > - (1-T)$ by \ben \label{defLTW}
L^{G}_T(\theta) := \int_0^T g^{G}(t, \theta)\!\ dt\,. \een
 If $V_\ell^*$ is equipped with the topology
$\sigma(V_\ell^*, V_\ell)$, the set \[B_a := \{\mu \in V_\ell^* :
|\mu|_{[0,T]} \leq a\}\] is compact according to the Banach-Alaoglu
theorem. Now $- \dot\Theta^G_n $ is a positive measure and its total
mass is $-\dot\Theta^G_n([0,T]) = -\Xi_n (T)$. We have then
\[\mathbb{P}\Big(\dot\Theta^G_n \notin B_a\Big) = \mathbb{P}\Big(\Theta^G_n (T) <
-a\Big)\,.\] Now for $\theta <0$
\[\mathbb{P}\left(\Theta^G_n (T) <
-a\right) \leq e^{\beta'\theta n^2 a}\!\ \mathbb{E}\exp
\{n^2\beta'\theta\Theta^G_n (T) \}\] so that, taking logarithm and
applying Lemma \ref{cvcgfG} we get, for $\theta \in (-(1-T),0)$
\[\limsup_n \frac{1}{\beta' n^2} \log \mathbb{P}\left(\Theta^G_n (T) <
-a\right) \leq \theta a + L_T^G(\theta)\,.\] It remains to let $a
\rightarrow \infty$ and we have proved the exponential tightness.

Let us note that the restriction  $T<1$ is crucial in the above
proof.

\subsubsection{Exposed points}
\label{exposed} Let ${\mathcal R}$ be the set of functions from
$[0,T]$ into $\mathbb {R}$ which are positive, continuous and with bounded
variation. Let ${\mathcal F}$ be the set of those $m \in V_\ell^*$
(identified with ${\mathcal M}_T$ as in \cite{leo1}) which are
absolutely continuous and whose density $\rho$ is such that $-\rho
\in {\mathcal R}.$ Let us prove that such a $m$ is exposed, with
exposing hyperplane $f_m(t) = \lambda(t, \rho(t))$ (recall
(\ref{lambdaH})). Actually we follow the method of  \cite{grz}. For
fixed $t$, $g^\star(t, .)$ is strictly convex on $(-\infty, 0)$ so
that, if $z\not= \rho(t)$, we have
$$g^\star(t, \rho(t)) - g^\star(t, z) < \lambda(t,\rho (t)) (\rho (t) -z)\,. $$
Let
 $d\xi = \tilde l(t) dt + \xi^\perp$  the Lebesgue decomposition of some element  $\xi \in {\mathrm M}_T$ such that
 $\widetilde I^G_{[0,T]}(\xi) < \infty$.
Taking $z = \tilde l(t)$ and integrating, we get \be \int_0^T
g^\star (t, \rho(t) dt - \int_0^T g^\star (t, \tilde l(t)) dt <
\int_0^T \lambda(t,\rho (t)) \rho (t) dt - \int_0^T \lambda(t,\rho
(t)) \tilde l(t) \ dt \ee and since $\int_0^T g^\star (t, \tilde
l(t)) dt = \int_0^T L_a^G (t, \tilde l(t)) dt \leq  \widetilde
I_{[0,T]}^G (\xi)$ this yields \be \widetilde I_{[0,T]}^G (m) -
\widetilde I_{[0,T]}^G (\xi)  < \int_0^T f_m dm -\int_0^T f_m
d\xi\,. \ee The following lemma says that this set of exposed points
is rich enough.
\begin{lem}
\label{exposition} Let $m \in V_r$ such that $\widetilde I_{[0,T]}^G
(m) < \infty$. There exists a sequence of functions $l_n \in
{\mathcal R}$
 such that
\begin{enumerate}
\item $\lim_n l_n(t) dt =  -m$ in $V_\ell^*$ with the $\sigma(V_\ell^*, V_\ell)$
topology,
\item  $\lim_n \widetilde I_{[0,T]}^G (-l_n(t) dt) = \widetilde I_{[0,T]}^G(m)\,.$
\end{enumerate}
\end{lem}
\prf The method may be found in \cite{grz} and in \cite{gamb}. The
only difference is in the topology because we want to recover
marginals. We will use  the basic inequality which holds for every
$\epsilon \leq 0$ : \ben \label{ineqH} L_a^G (t, v+ \epsilon) \leq
L_a^G (t,v) - \epsilon(1-t) \een Let $m= m_a + m_s$ such that
$\widetilde I_{[0,T]}^G (m) < \infty$. From (\ref{42}) and
(\ref{LaLs}) it is clear that $-m_a$ and $-m_s$ must be  positive
measures.
\medskip

\noi{First step} We assume that $m = -l(t) dt - \eta$ with $l \in
L^1([0,T]; dt)$ and $\eta$ a singular positive measure.
 One can find a sequence of  non negative continuous functions
 $h_n$  such that $h_n(t) dt \rightarrow \eta$
for the topology $\sigma(V_\ell^*, V_\ell)$. Indeed every function
$\psi \in V_\ell$ may be written as a difference
 $\psi_1 - \psi_2$ of two increasing functions. There exists a unique (positive) measure $\alpha_{1}$
 such that $\psi_1(t) = \alpha_{1} ([t,T])$ for every $t \in [0,T]$. Moreover,
 the function $g = \eta([0,\cdot]) \in V_r$ is non decreasing and may be approached by a sequence
  of continuously derivable and non decreasing functions $(g_n)$ such that
  $g_n \leq g$. Setting $h_n := g'_n$ and
 $\nu_n = h_n (t) dt$, the dominated convergence theorem gives
\be \langle\psi_{_1}, \nu_n\rangle = \int_0^T \nu_n
([0,t])\alpha_{1} (dt) \rightarrow \int_0^T \eta([0,t]) \alpha_{1}
(dt) \,. \ee With the same result for  $\psi_2$ we get \be
\langle\psi, \nu_n\rangle &=& \int_0^T \nu_n ([0,t])\alpha_{1} (dt)
- \int_0^T \nu_n ([0,t])\alpha_{2} (dt)
\\&\rightarrow& \int_0^T \eta([0,t]) \alpha_{1} (dt) - \int_0^T
\eta([0,t]) \alpha_{2} (dt) \,. \ee or $\lim_n \langle\psi,
\nu_n\rangle = \langle\psi, \eta\rangle$. On the one hand, the lower
semi-continuity of
  $\widetilde I_{[0,T]}^G$ yields
\be \liminf_n \widetilde I_{[0,T]}^G \left(-(l(t) +
h_n(t))dt\right)\geq \widetilde I_{[0,T]}^G (m)\,. \ee On the other
hand, integrating (\ref{ineqH}) yields \be
\widetilde I_{[0,T]}^G (-(l(t) + h_n(t))dt)&\leq& \int_0^T L_a^G (t, -l(t)) dt  + \int_0^T (1-t)h_n (t) dt\\
\nonumber &\rightarrow& \int_0^T L_a^G (t, -l(t)) dt  + \int_0^T
(1-t) \eta (dt) = \widetilde I_{[0,T]}^G (m)\,. \ee
\medskip

\noi{Second step} Let us assume that $m = -l(t) dt$ with $l \in
L^1([0,T]; dt)$ and for every $n$, let us set $l_n = \max(l, 1/n)$.
It is clear that as $n \rightarrow \infty$, then $l_n \downarrow l$.
On the one hand the lower semi-continuity gives \[\liminf_A
\widetilde I^G_{[0,T]} (-l_n(t)dt) \geq I_{[0,T]}^G (-l(t)dt)\,.\]
On the other hand, by integration of inequality (\ref{ineqH}), since
$l_n -l \leq 1/n$
$$I^G_{[0,T]} (-l_n(t)dt) \leq I^G_{[0,T]} (-l(t)dt) + \frac{1}{n}\,.$$
It is then possible to reduce the problem to the case of functions
bounded below.
\medskip

\noi{Third step} Let us assume that $m = -l(t) dt$ with $l \in
L^1([0,T]; dt)$ and bounded below by $A >0$. One can find a sequence
of continuous  functions $(h_n)$ with bounded variation such that
$h_n \geq A/2$
 for every $n$ and such that
 $h_n \rightarrow l$ a.e. and in $L^1 ([0,T] , dt)$.
We have $h_n(t) dt \rightarrow l(t)dt$ in $\sigma(V_\ell^*, V_\ell)$
and since $L_a^G (t, \cdot)$ is uniformly Lipschitz on $(-\infty,
-A/2]$, say with constant $\kappa$, we get \be |\widetilde
I_{[0,T]}^G (-h_n(t) dt) - \widetilde I_{[0,T]}^G (-l(t)dt)|\leq
\kappa \int_0^T |h_n (t) - l(t)| dt \rightarrow 0\,. \ee Actually,
$h_n \in {\mathcal R}$ and $\varphi_n (t) := \lambda(t, -h_n(t))$
satisfies the assumption of Lemma \ref{cvcgfG} since
$$1 +  \varphi_n (t) -t\geq \frac{t}{1-e^{-A/2}}\,.$$

\subsubsection{End of the proof of Theorem \ref{LDPH}}
\label{end}  The first step is the upper bound for compact sets. We
use Theorem 4.5.3 b) in \cite{DZ} and the following lemma.
\begin{lem}
\label{petitlemme} For every $\delta > 0$ and $m \in V_\ell^*$,
there exists $\varphi_\delta$ fulfilling the conditions of Lemma
\ref{cvcgfG} and such that \ben \label{condZani} \int_0^T
\varphi_\delta dm - \Lambda_T^G (\varphi_\delta) \geq \min \Big[
I_{[0,T]}^G (m) - \delta ,\  \delta^{-1}\Big]\,. \een
\end{lem}
The second step is the  upper bound for closed sets : we use the
exponential tightness. The third step is the  lower bound for open
sets. The method is classical (see \cite{DZ} Theorem 4.5.20 c)),
owing to Lemma \ref{exposition}.
\medskip

To prove Lemma \ref{petitlemme}, we start from the definition
(\ref{dualeG1}) or (\ref{dualeG2}). One can find $\bar
\varphi_\delta \in V_\ell$ satisfying (\ref{condZani}). If $\bar
\varphi_\delta$ does not check assumptions of the lemma we add
$\varepsilon > 0$ to $\bar \varphi_\delta$ which allows to check
them and satisfy
 (\ref{condZani}) up to a change of $\delta$. $\Box$

\subsection{Proof of Theorem \ref{margH}}

We use the contraction from the LDP for paths. Since the mapping $m
\mapsto m([0,T])$ is continuous from $D$ to $\mathbb{R}$, the family
$\Theta^G_n(T)$
 satisfies the
LDP with good rate function specified by (\ref{opt1G}):
$$I^{G}_T (\xi) = \inf \{ I_{[0,T]}^G (v) \ ; \  v(T) = \xi\}\,.$$
Since the process $\Theta^G_n $ takes its values in $(-\infty, 0]$
(remember Hadamard inequality), it is clear that  $I^{G}_T (\xi) =
\infty$ for $\xi > 0$.

Fixing $\xi < 0$, we can look for optimal $v$. Let $\theta  > -
(1-T)$ (playing the role of  a Lagrange multiplier).

By the duality property (\ref{identifH})
\[
g^\star \Big(t , \frac{d\dot v_a}{dt}(t)\Big)\geq \theta \frac{d\dot
v_a}{dt}(t) - g^G(t, \theta) \,.
\]
 Integrating and using (\ref{tauxHTG}), (\ref{gLa}) and (\ref{98G}) we get 
\ben I_{[0,T]}^G (v) \geq \theta \dot v_a ([0,T]) - L^G_T (\theta)-
\int_0^T (1-t) \ d\dot v_s (t)\,, \een For every $v$ such that $v(T)
= \xi$ it turns out that \ben \label{riche} I_{[0,T]}^G (v) \geq
\theta \xi -L^G_T (\theta) - \int_0^T (1-t + \theta) \ d\dot v_s
(t)\geq \theta \xi -L^G_T (\theta) \,. \een Besides, from
(\ref{lambdaH}) the ordinary differential equation \be
\lambda^G(t, \phi'(t))&=&\theta\\
\phi(0) &=& 0\,, \ee admits for unique solution in ${\mathcal
C}^1([0,T])$
\[
t \mapsto \phi (\theta; t) := {\mathcal J}(1+\theta)-{\mathcal
J}(1-t +\theta)  -t \log (1+\theta)\,.
\]
Now, since
$$\frac{\partial}{\partial \theta} \phi (\theta; T) =
-\left[\log\Big(1 - \frac{T}{1+\theta}\Big) +
\frac{T}{1+\theta}\right] > 0$$ we see that the mapping $\theta
\mapsto \phi(\theta; T)$ is bijective from $[-(1-T), \infty)$ onto
$[-T, 0)$. Moreover, by duality
$$g^\star\Big(t , \frac{\partial}{\partial t}\phi(\theta,t)\Big) =
\theta\frac{\partial}{\partial t}\phi(\theta,t) -g^G(t, \theta)\,.$$
There are two cases.

$\bullet$ If $\xi \in [- T, 0)$, there exists a unique  $\theta_\xi$
such that $\phi(\theta_\xi, T) = \xi$ (i.e. the relation
(\ref{=xi1H}) is satisfied). For $v^\xi := \phi(\theta_\xi \!\ ,
\cdot)$, we get from (\ref{tauxHTG}), (\ref{gLa}) and (\ref{defLTW})
again \be I_{[0,T]}^G(v^\xi) &=& \theta_\xi \xi - L^G_T(\theta_\xi)
\ee so that $v^\xi$ realizes the infimum in (\ref{opt1G}). A simple
computation ends the proof of the first statement of Theorem
\ref{margH}.

 Let us note that at the end point $\xi = -T$, we have
 \be
 \theta_\xi = -(1-T) \  , \ v^\xi (t)= {\mathcal J}(T) -{\mathcal J}(T-t)  - t \log T\ , \
 (v^\xi)'(t) = \log (1-t/T)\,.
 \ee
 Finally
\be
I_T^{G}(-T) &=& 2T(1-T) + \left(  F(1) - F(1-T)-F(T) + T^2 \log T\right)\\
\nonumber &=& \frac{T(1-T)}{2} + \frac{T^2 \log T}{2} -
\frac{(1-T)^2 \log (1-T)}{2} + \frac{3}{4} \,. \ee

$\bullet$ Let us assume $\xi= -T - \varepsilon$ with $\varepsilon >
0$. Plugging $\theta = - (1-T)$ in (\ref{riche}) yields, for every
$v$ such that $v(T) = \xi$ \be I_{[0,T]}^G (v) \geq - (1-T)\xi
-L^G_T \left(-(1-T)\right) =   (1-T) \varepsilon + I_T^{G}(-T)
\,,\ee and this lower bound is achieved by the measure $\widetilde v
= (v^{-T})'(t) dt - \varepsilon \delta_T(t)$, since
\[\int_0^T L_a^G \left(t, (v^{-T})'(t)\right) dt = I_T^{G}(T)\ \ ,
\ \ \int_0^T  (1-t)  \!\ \varepsilon \!\ d\delta_T(t) = (1-T)\!\
\varepsilon\,.
\]
It remains to look at $\xi = 0$. Taking  $\xi = 0$ in (\ref{riche}),
we get
$$I_T^G (0) \geq -L^G_T (\theta)$$
for every $\theta \geq -(1-T)$. Now, from (\ref{defg}) and
(\ref{defLTW}) 
 we may write
\be-L^G_T (\theta) = \int_0^T (1-t) \log (1-t)\ dt + \int_0^T
(1+\theta)\log\Big(1- \frac{t}{1+\theta}\Big)\!\ dt\\ + \int_0^T t
\log (1-t + \theta)\!\ dt\,.\ee When $\theta$ tends to infinity, the
second term tends to zero and the third, which is bounded above by
$(T^2/2)\log(1-T+\theta)$
tends to infinity. Finally $I_T^G (0) = \infty$.

That ends the proof of the second statement of Theorem \ref{margH}.
$\Box$

\begin{remark}
\label{T=1} It is possible to try  a direct method to get
(\ref{idetH}), (\ref{mauvH}) using  G\"artner-Ellis' theorem
(\cite{DZ}, Theorem 2.3.6). From Lemma \ref{cvcgfG} the limiting
n.c.g.f. of $\Theta^G_n(T)$ is $L^G_T$ which is analytic on $(-
(1-T), \infty)$. When $\theta \downarrow -(1-T)$, we have
$(L^G_T)'(\theta) \downarrow -T$. We meet  a case of so called non
steepness. To proceed in that direction we could use the method of
 time dependent change of probability
(see  \cite{DZT}). We will not give details here.
Nevertheless, this approach allows to get one-sided large deviations
in the critical case $T=1$. Actually we get \be \lim_n
\frac{1}{\beta'n^2} \log \mathbb{P}( \Theta^G_{n}(1) \geq nx) = - I_1^G (x)
\ee for $x \geq -1$. The value $x=-1$ corresponds to the limit of
$\Theta^G_{n}(1)$. note that the second (right) derivative of
$I_1^G$ at this point is zero (or equivalently $\lim (L^G_1){''}
(\theta) = \infty$ as $\theta \downarrow 0$) , which is consistent
with previous results on variance. I do not know the rate of
convergence to $0$ of   $\mathbb{P}( \Theta^G_{n}(1) \leq nx)$ for $x< -1$.
\end{remark}

\subsection{Proof of Theorem \ref{LDPwish} and Theorem \ref{margW}}

\label{similarroute} Again, the three routes are possible to tackle
the problem of large deviations for determinant of Wishart matrices.
A direct method would use the cumulant generating function from
(\ref{newcgf}) and would meet computations similar to those seen in
the Uniform Gram case.

 To avoid repetitions, we use the decomposition (\ref{relproc}), drawing benefit from
an auxiliary study of  $S_{n,r}$.
\begin{lem}
\label{LDPaux} The sequence  $\{ n^{-1}\!\ S_{n}(t) , t \in [0,
1)\}_n$ satisfies a LDP in the space $(D, \sigma(D, {\mathrm
M}_<))$ in the scale $2\beta^{-1}n^{-2}$ with good rate function
\ben \label{59} I_{[0,1)}^{S}(v) = \int_{[0,1)}
L_a^{S}\Big(\frac{d\dot v_a}{dt} (t)\Big) dt + \int_{[0,1)} L_s^{S}
\Big(\frac{d\dot v_s}{d\mu} (t) \Big) d\mu(t) \een where \ben
\label{60} L_a^{S} (y) = \left(e^{y} -y-1\right) \ , \ L_s^S (y) = -
y \delta(y|(-\infty,0))\,, \een and  $\mu$ is any  measure
dominating $d\dot v_s$.
\end{lem}
Let us stress that the instantaneous rate functions  are time
homogeneous and then we may write $[0,1]$ instead of $[0,1)$.
\subsubsection{Proof of Lemma \ref{LDPaux}}
It is a route similar to the proof of Theorem \ref{LDPH} in Section
\ref{7.4G} (see also  \cite{Najim1}). We start from
(\ref{defsigmaproc}) so that
\[\frac{1}{n} \dot S_n = \sum_{j=1}^n \Big(\log\varepsilon_k^{(n)}\Big)\!\ \delta_{j/n}\,.\]
Withe help of (\ref{newcgf}) this yields : \be \log \mathbb{E}  \exp< \beta'
n \dot S_n , \gamma
>
  = \sum_{k=1}^n \left[\ell\left(\beta'n \Big(1+\gamma\Big(\frac{k}{n}\Big)\Big)\right)  -
  \beta'\gamma\Big(\frac{k}{n}\Big) \log(\beta'n)\right]\\
  - n \ell(\beta'n)
  \ee
if $\gamma(s) + 1 > 0$ for every $s \in [0,1]$. A little computation
shows that the limiting n.c.g.f. is \ben \label{128} {\mathcal
L}^{S} (\gamma) = \int_0^1 {\mathcal J}(1+\gamma(t)) dt\,, \een
which yields (\ref{60}) by duality (see \cite{Rocky1} again). $\Box$
\subsubsection{Proof of Theorem  \ref{LDPwish}}
Let  $\Theta^L_n = n^{-1} \log \Delta_n^{L, \beta}$. We deduce from
Lemma \ref{LDPaux} and Theorem \ref{LDPH} that the sum $\dot
\Theta^L_n = \dot \Theta^G_n + \frac{1}{n}\dot S_n$ satisfies a LDP
in the same scale with good rate function $I_{[0,T]}^{G}\Box
I_{[0,T]}^{S}$ where $\Box$ denotes the infimum convolution :
\[(f\Box g)(x) = \inf\{f(x_1) + g(x_2)\ |\ x_1 + x_2 = x \} \,.\]
The two characteristics of the rate function are then
\begin{eqnarray*}
L_a^L
 &=& \inf_v \{L_a ^{G}(v) + L_a^{S}(u-v)\}\\
L_s^L &=& \inf_v \{L_s^{G}(v) + L_s^{S}(u-v)\}\,.
\end{eqnarray*}
which yield (\ref{identif}) by an explicit computation. $\Box$
\bigskip

Alternatively, it is possible to sum the two n.c.g.f. ((\ref{98G})
and (\ref{128})) and get the rate function by  duality.  We
 claim : if $\gamma(s) +1
> 0$ for every $s \in [0,1]$ \ben \label{ncgfL} \frac{1}{\beta'n^2}
\log \mathbb{E}\exp \langle\beta' n^2\dot\Theta_n^L , \gamma\rangle
\rightarrow \int_0^T g^L(t, \gamma(t))\!\ dt\,, \een where \ben
g^L(t, \gamma) = g^G(t, \gamma) + {\mathcal J}(1 + \gamma) =
{\mathcal J}(1 -t + \gamma) - {\mathcal J}(1 -t)\,. \een

\subsubsection{Proof of Theorem \ref{margW}}

We may either use the contraction $\Theta_n^L \mapsto \Theta_n^L (T)$ 
 or
establish a LDP for the marginal $S_n(T)$ and then perform an
inf-convolution. We leave the details of the proof to the reader. We
just give the expression of the optimal path when it exists.

For $\theta > -(1-T)$, the function
\[
t \mapsto \phi (\theta; t) := {\mathcal J}(1+\theta)-{\mathcal
J}(1-t +\theta)  \,.
\]
is in ${\mathcal C}^1([0,T])$ and the mapping $\theta  \mapsto
\phi(\theta; T)$ is bijective from $[-(1-T), \infty)$ onto $[\xi_T,
\infty)$, where $\xi_T = {\mathcal J}(T) - 1$.

Fixing $\xi\geq\xi_T$, we can look for optimal $v$. There exists a
unique $\theta_\xi$ such that $\phi(\theta_\xi, T) = \xi$. Then
$v^\xi := \phi(\theta_\xi \!\ , \cdot)$ is the optimal path ($v^\xi$
realizes the infimum in (\ref{opt1L}). Let us note that at the end
point $\xi = \xi_T$, we have
 \be
 \theta_\xi = -(1-T) \  , \ v^\xi (t)= {\mathcal J}(T) -{\mathcal J}(T-t)  \ , \ (v^\xi)'(t) = \log (T-t)\,. \Box
 \ee

\begin{remark}
\label{T=1L} It is possible to get (\ref{idetW}), (\ref{mauvW})
using  G\"artner-Ellis' theorem (\cite{DZ}, Theorem 2.3.6). We are
in the same situation as in Remark \ref{T=1}.

This approach allows to get one-sided large deviations in the
critical case $T=1$. Actually we get \be \lim_n \frac{2}{\beta n^2}
\log \mathbb{P}( \log \Delta^{L, \beta}_n (1) \geq n x) = - I_1^L (x) \ee
for $x \geq -1$. The value $x=-1$  corresponds to the limit.
 Note that the second (right) derivative of $I_1^L$ at
this point is zero (or equivalently $\lim (L^L_1)'' (\theta) =
\infty$ as $\theta \downarrow 0$), which is consistent with
previous results on variance. We do not know the rate of convergence
to $0$ of   $\mathbb{P}( \log \Delta^{L, \beta}_n (1) \leq nx)$ for $x< -1$.
\end{remark}
\subsection{Proof of Theorem \ref{LDPpath} and Theorem \ref{margJ}}
\label{7.3} We may try again to use the beta-gamma algebra, but  we
do not succeed to go until the end. Let as in Subsection \ref{55},
$U_n$ and $V_n$ be the two Laguerre variables. From the exponential
tightness of $U_n$ and $V_n$, we deduce easily the exponential
tightness of $W_n$. From \cite{MR1851048}, the sequence
$W_n$ contains subsequences satisfying LDP. If for such a
subsequence we call $I^p$ the rate function, the independence gives
\[I^U = I^V \Box I^p\]
This equation has many solutions and only one convex solution, which is 
 \be I^p = I^U\boxminus I^V
\ee defined by
\[(f\boxminus g)(x) = \sup\{f(x_1) - g(x_2) \ | \ x_1-x_2 = x\}\]
( \cite{Mazure}). But we do not know a priori that $I^p$ is convex.

  We choose to use the beta-gamma trick to study
the n.c.g.f. For the remaining we do not give details since it is
similar to the above cases  and again based on the ideas of Baldi
theorem (\cite{DZ}) and a variational formula.

\subsubsection{Convergence of the n.c.g.f.}
Put $\Theta_n^J = n^{-1} \log \Delta_n^{J, \beta}$ so that
\be \dot\Theta_n^J = \frac{1}{n} \sum_{k=1}^{n_1}
\Big(\log \rho_{j,n}^{J,\beta}\Big)\!\ \delta_{j/n}\,,
 \ee
 and put for $T \leq \tau_1$ and $\varphi\in V_\ell^{T}$ :
\be
{\mathcal
L}^J_{n, \lfloor nT\rfloor}(\varphi) = \frac{2}{\beta n^2} \log
\mathbb{E}\exp
n\langle\dot\Theta^J_n , \varphi\rangle 
 \,.
\ee
\begin{lem}
\label{cvcgfH} If $\varphi \in V_\ell^{\tau_1}$ satisfies
$\varphi(s) +\tau_1 -s > 0$ for every $s\in (0,T]$, then \ben
\label{98J} \lim_n {\mathcal L}_{n, \lfloor nT\rfloor}^J (\varphi)
=\Lambda_{[0,T]}^J (\varphi) := \int_0^T g^J(s, \varphi(s)) \ ds \,,
\een where, for $\theta +\tau_1 -s > 0$ \ben \label{defgJ} g^J(s,
\theta) ={\mathcal E}\left(\tau_1-s+\theta, \tau_2, \theta\right)
\,. \een
\end{lem}
\prf  From (\ref{ThetamoinsTheta}) we have \be \langle n
\dot\Theta_n^J, \gamma \rangle + \langle (n_1 +
n_2) \dot\Theta^L_{n_1 + n_2}, \gamma((n_1+n_2)\cdot/n)\rangle=\\
\langle n_1\dot\Theta^L_{n_1}, \gamma(n_1\cdot/n)\rangle +
\log\frac{n_1}{n_1+n_2} \sum_{k=1}^{\lfloor nT\rfloor}
\gamma(k/n)\ee and then, by independence, \be \log \mathbb{E} \exp \langle
\beta' n^2 \dot\Theta^{J}_n , \gamma\rangle = \log \mathbb{E} \exp \langle
\beta' nn_1 \dot\Theta^L_{n_1}, \gamma (n_1\cdot/n)\rangle\ \ \ \ \
\ \ \ \ \ \ \ \ \ \ \ \ \ \ \ \ \
\\-
   \log \mathbb{E} \exp
\langle \beta' n(n_1 + n_2) \dot\Theta^L_{n_1+n_2} ,
\gamma((n_1+n_2)\cdot/n)\rangle \\+ n \log\frac{n_1}{n_1+n_2}
\Big(\sum_{k=1}^{\lfloor nT\rfloor} \gamma(k/n)\Big)\,.\ee By a
slight modification  of (\ref{ncgfL}) we have, for $p/n \rightarrow
\tau$ \ben \label{scale} \frac{1}{\beta' p^2} \log \mathbb{E} \exp \langle
\beta' np \dot\Theta^L_r , \gamma(p\cdot/n)\rangle \rightarrow
\frac{1}{\tau}\int_0^T g^L\Big(\frac{s}{\tau},
\frac{\gamma(s)}{\tau}\Big) \ ds\,, \een so that taking
$\tau=\tau_1$ and $\tau = \tau_1 + \tau_2$ and subtracting, we get
\be \frac{1}{\beta' n^2} \log \mathbb{E} \exp \langle \beta' n^2
\dot\Theta^J_n , \varphi\rangle \rightarrow \int_0^ {T} g^J\left(s,
\gamma(s)\right)\,,\ee where \be g^J\left(s, \theta\right)= \tau_1
g^L\Big(\frac{s}{\tau_1}, \frac{\theta}{\tau_1}\Big) - (\tau_1+
\tau_2) g^L\Big(\frac{s}{\tau_1+\tau_2},
\frac{\theta}{(\tau_1+\tau_2)}\Big)
\\+ \theta\!\ \log\frac{\tau_1}{\tau_1+\tau_2}\,,\ee and this is
equivalent to (\ref{defgJ}). $\Box$
\subsubsection{Duality}
Let us define $\Lambda_{[0,T]}^J (\varphi) = +\infty$ if $\varphi$
does not satisfy the assumption of Lemma \ref{cvcgfH}. The dual of
$\Lambda_{[0,T]}^J$ is then \ben \label{dualeH1}
\Big(\Lambda_{[0,T]}^J\Big)^\star (\nu) = \sup_{\varphi \in V_\ell}
\left\{<\nu, \varphi> - \Lambda_{[0,T]}^J (\varphi) \right\} \een
for $\nu \in V_\ell^*$. Mimicking the method of 
\cite{leo1} p. 112-113, we get \ben \label{dualeH2}
\Big(\Lambda_{[0,T]}^J\Big)^\star (\nu) = \sup_{\varphi \in
{\mathcal C}} \left\{<\nu, \varphi> - \Lambda_{[0,T]}^J (\varphi)
\right\} \een where ${\mathcal C}$ is the set of continuous
functions from $[0,T]$ into $\mathbb {R}$ vanishing at $0$. Then we apply
Theorem 5 of  \cite{Rocky1}. We get \ben \label{rocky}
\Big(\Lambda_{[0,T]}^J\Big)^\star (\nu) = \int_0^T
\Big(g^J\Big)^\star\Big(t, \frac{d\nu_a}{dt}\Big) \!\ dt + \int_0^T
r^J\Big(t, \frac{d\nu_s}{d\mu}(t)\Big)\!\ d\mu(t) \een where \ben
\label{defgJstar} \Big(g^J\Big)^\star (s, y) = \sup_\lambda
\left\{\lambda y - g^J(s,\lambda)\delta(\lambda | (-\tau_1 ,
\infty)) \right\}\,. \een This supremum is achieved by \ben
\label{lambdaJ} \lambda^J(s,y) = -(\tau_1 -s)
 + \frac{\tau_2
 }{e^{-y}-1}
\een
and we have
\ben
\label{dual}
\Big(g^J\Big)^\star (s, y) &=& \lambda^J(s,y) y - g^J(s, \lambda^J(s,y))\\
\label{59J} &=& \left(\tau_1 + \tau_2
 -s\right){\bf H}\Big(\frac{\tau_1 -s
  }{\tau_1 + \tau_2
   -s}\ \Big|\ e^y\Big)\,.
\een
The recession is $r^J(s;y) = -(\tau_1 -s)y$ if $y < 0$.

\subsubsection{Proof of Theorem \ref{margH}}
 We use the contraction from the LDP for paths.
Since the mapping $m \mapsto m([0,T])$ is continuous from $D$ to
$\mathbb{R}$, the family $\{\Theta^J_{n}(T)\}_n$ satisfies the   LDP with
good rate function given by (\ref{opt1J}). Since the process
$\Theta_n $ takes its values in $(-\infty, 0]$ , it is clear that
$I^{J}_T (\xi) = \infty$ for $\xi > 0$. Fixing $\xi < 0$, we can
look for optimal $v$, i.e. a path $(v(t), t \in [0,T])$ such that
$v(T) = \xi$ and $v$ achieves the infimum in (\ref{opt1J}). Fix
$\theta \geq T-\tau_1$ (playing the role of  a Lagrange multiplier).
 In view of (\ref{rocky}), (\ref{defgJstar}) and (\ref{lambdaJ}), it is clear that (in the generic case) the Euler-Lagrange equation is
\be
\lambda^J(s, \phi'(s))&=&\theta\\
 \phi(0) &=& 0\,. \ee This ordinary differential equation admits
for unique solution in ${\mathcal C}^1([0,T])$
\[
s \mapsto \phi^J (\theta; s) := {\mathcal E}(\theta+\tau_1, \tau_2,
s)\,.
\]
To know if the path $\phi^J$ may have $\xi$ as its terminal value,
look at \be{\mathcal E}'(\theta + \tau_1, \tau_2, T) &=&
\frac{\partial}{\partial \tau} {\mathcal E} (\tau, \tau_2,
T)|_{\tau= \theta+\tau_1}\\ &=& \log \Big(1 - \frac{T}{\theta +
\tau_1+\tau_2}\Big) - \log \Big(1 - \frac{T}{\theta +
\tau_1}\Big)\,; \ee since it is positive, we see that the mapping
$$\theta \longmapsto {\mathcal E}(\theta + \tau_1 , \tau_2 ,T)
$$
is continuous and increasing from $[T-\tau_1, \infty)$ onto
${\mathcal D}_T= [\xi_T^J , \ 0)$. If $\xi \in [\xi_T^J, 0)$, we
call $\theta_\xi$ the unique solution of   $\phi^J(\theta, T) = \xi$
or in other words, \be {\mathcal E}(\theta_\xi + \tau_1, \tau_2 ,
T)= \xi\,, \ee and we set $v^\xi := \phi^J(\theta_\xi \!\ , \cdot)$.

To end the proof, let us now consider some inequalities.  The
duality property (\ref{defgJstar}) gives, for every $v$ and $t$ \ben
\label{tointegrate} \Big(g^J\Big)^\star \Big(t , \frac{d\dot
v_a}{dt}(t)\Big)\geq \theta \frac{d\dot v_a}{dt}(t) - g^J(t, \theta)
\,. \een Setting \ben \label{defLT} L_T^J (\theta) := \int_0^T g^J
(t, \theta)\ dt\,, \een integrating (\ref{tointegrate}) and using
(\ref{intact0}), (\ref{LaJ}) and (\ref{defLT}) we get \be
I_{[0,T]}^J (v) \geq \theta \dot v_a ([0,T]) - L^J_T (\theta)-
\int_0^T (\tau_1-t) \ d\dot v_s (t)\,. \ee For every $v$ such that
$v(T) = \xi$ it turns out that \ben \label{richeJ} I_{[0,T]}^J (v)
\geq \theta \xi -L_T^J (\theta) - \int_0^T (\tau_1-T+  \theta) \
d\dot v_s (t)\geq \theta \xi -L_T^J (\theta) \,. \een There are
three cases.

$\bullet$ If $\xi \in
[\xi_T^J, 0)$
, we get
\be
I_{[0,T]}^J(v^\xi) &=& \theta_\xi \xi - L_T^J(\theta_\xi)
\ee
so that $v^\xi$ realizes the infimum in (\ref{opt1J}). A simple computation leads to (\ref{idetH}) which
ends the proof of the first statement of Theorem \ref{margH}.

 Let us note that at the end point $\xi = \xi_T^J$, we have
 \be
 \theta_\xi =   (T-\tau_1)\  , \ v^\xi (t)= {\mathcal E}(T, \tau_2, t)\ , \ (v^\xi)'(t) = \log \frac{T-t}{\tau_1 + \tau_2 -t}\,.
 \ee

$\bullet$ Let us assume $\xi=  \xi_T^J -\varepsilon$ with
$\varepsilon > 0$. Plugging $\theta = T-\tau_1$ in (\ref{richeJ})
yields, for every $v$ such that $v(T) = \xi$ \be I_{[0,T]}^J (v)
\geq (T-\tau_1)\xi_T^J -L_T^J(T-\tau_1) -\varepsilon(T-\tau_1) =
I_T^J (\xi_T^J)   -\varepsilon (T-\tau_1)\,, \ee and this lower
bound is achieved by the measure  $\tilde v =
\Big(v^{\xi_T^J}\Big)'(t)dt -\varepsilon d\delta_T(t)$, since
$$\int_0^T L_a^J\Big(t,\Big(v^{\xi_T^J}\Big)'(t)\Big)dt = I_T^J(\xi_T^J)\ \ , \ \ \int_0^T (\tau_1 -t)\varepsilon d\delta_T(t) = (\tau_1 -T)\varepsilon\,.$$
$\bullet$ It remains to look at $\xi = 0$. Taking $\xi = 0$ in
(\ref{richeJ}), we get $I^J_{[0,T]}(0)\geq -L_T^J(\theta)$ for every
$\theta \geq T - \tau_1$. Now, from (\ref{defgJ}) and (\ref{defLT}),
we may write (after some calculation) \ben \nonumber -L_T^J(\theta)
= \int_0^T -{\mathcal E}(\tau_1 - t + \theta, \tau_2, \theta)dt &=&
\int_0^T \int_0^\theta \log \Big(1 + \frac{\tau_2}{\tau_1 -t
+s}\Big)ds\\\nonumber &\geq& T \int_0^\theta \log \Big(1 +
\frac{\tau_2}{\tau_1  +s}\Big)\ ds\,, \een which tends to infinity
as $\theta \rightarrow \infty$. We conclude $I^J_{[0,T]}(0) =
\infty$. $\Box$
\section{Proofs of Theorems of Section \ref{contrac}}
\subsection{Proof of Proposition  \ref{infIW}}Let $\theta \in \mathbb{R}$  be a Lagrangian factor. We begin  by minimizing
\be
 I^{spL}_T (\mu)- \theta T \int \log x \ d\mu(x)  = T^2
\left[- \Sigma(\mu) + 2 \int q_{\lambda,s} (x)
 \ d\mu(x)\right] + 2 B(T)
\ee where \ben \label{lambdas} q_{\lambda, s} (x) = \lambda x - s
\log x, \ \ \lambda = \frac{1}{2T} , \ s = \frac{1-T +
\theta}{2T}\,. \een In the book of \cite{Saff} p.43 example 5.4,
it is stated that  for $\lambda > 0$ and $2s+1 > 0$ fixed, the
infimum
$$\inf_\mu \ - \Sigma(\mu) + 2 \int q_{\lambda,s} (x) \ d\mu(x)$$
  is achieved by the unique
extremal measure
$\pi_{\sigma^2}^c$
 with
\be \sigma^2 = \frac{2s+1}{2\lambda} \ , \ c = \frac{1 }{2s+1}\,.\ee
We see from (\ref{lambdas}) that if $\theta > -1$ we can take: \be
\sigma^2 = 1 +\theta , \ \ \ c =
\frac{T}{\sigma^2}=\frac{T}{1+\theta}\,. \ee Now it remains to look
for $\theta$ such that the constraint $\mu \in {\mathcal A}(\xi/T)$
is saturated. Since \be \int \log x \ d\pi_{\sigma^2}^c (x) = \log
\sigma^2 + \int \log x \ d\pi_1^c(x) dx\,, \ee and thanks to
(\ref{often}) we see that $\theta$ must satisfy \be \xi = T \log
\sigma^2 - T \frac{{\mathcal J}(1-c)}{c} = {\mathcal J}(1+\theta) -
{\mathcal J}(1-T + \theta)\,, \ee which is exactly exactly
(\ref{=xi1W}).

To compute $I^{spL}_T(\pi_{\sigma^2}^c)$, we start from the
definition (\ref{HP1}): \be  I^{spL}_T(\pi_{\sigma^2}^c) = -
T^2\Sigma(\pi_{\sigma^2}^c) + T\int (x - (1-T)\log x)\
d\pi_{\sigma^2}^c(x) + 2 B(T)\,, \ee and transform
$\pi_{\sigma^2}^c$ to $\pi_1^c$ using the dilatation. In particular,
(\ref{HP2}) yields $$\Sigma(\pi_{\sigma^2}^c) = \log \sigma^2 +
\Sigma(\pi_1^c)$$ and $\Sigma(\pi_1^c)$ may be picked from formula
(13) p.10 in \cite{hiai1} :
$$\Sigma(\pi_1^c) = -1 + \frac{1}{2}\left(c^{-1} + \log c\ + (c^{-1} -1)^2 \log(1-c)\right)\,.$$
Besides   we have easily $\int x \ d\pi_1^c (x) = 1$. After some
tedious but elementary computations we get exactly the RHS of
(\ref{idetW}), which yields
\[I_T^{spL}(\pi_{\sigma^2}^c)= I_T^L(\xi)\,,\]
and ends the proof of (\ref{defxiL}). $\Box$
\subsection{Proof of proposition  \ref{infI}}
Let $\theta < t - \tau_1$  (Lagrangian multiplier). We begin  by
minimizing \ben\label{lastmin} I_t^{spJ}(\mu)\! &-&\!  \theta t \int \log x \
d\mu(x)
\\\nonumber&=& t^2 \left[-\Sigma(\mu) - 2\zeta_1 \int \log x\!\ d\mu(x) - 2
\zeta_2 \int \log(1-x)\!\  d\mu(x) \right] + C \een where $$2\zeta_1
=\frac{\tau_1 + \theta -t}{t}\  , \ 2\zeta_2 =\frac{\tau_2 -t}{t}\ \
\hbox{and}\  C = t^2 B\Big(\frac{\tau_1 -t}{t}, \frac{\tau_2
-t}{t}\Big)\,.$$
We use the following lemma.
\begin{lem}
\label{min} For $\zeta_1, \zeta_2 > 0$, the infimum of
$$ -\Sigma(\mu) - 2\zeta_1\int \log x \!\ d\mu(x)-2\zeta_2\int \log (1-x) \!\ d\mu(x)$$
among the probability measures $\mu$ on $[0,1]$ is achieved by
$\pi_{a_-, a_+}$ where \be (a_-,a_+) = \lambda_\pm (s_-, s_+) \ee
with
\[s_- = \frac{1 + 2\zeta_1}{2(1+\zeta_1 + \zeta_2)}\ , \
s_+ = \frac{1 + 2\zeta_1 + 2\zeta_2}{2(1+\zeta_1 + \zeta_2)}\,.
\]
\end{lem}
The infimum in (\ref{lastmin}) is achieved by
$\pi_{\tilde\xi, \tilde\eta}$, where
\[
(\tilde\xi, \tilde\eta) = \lambda_\pm(\tilde s_- , \tilde s_+)\,, \  \ \tilde s_- = \frac{\tau_1 + \theta}{\tau_1 + \tau_2 +\theta}\ , \ \tilde s_+ = \frac{\tau_1 + \theta + \tau_2  -t}{\tau_1 + \tau_2 +\theta}\,.
\]
It should be clear that \be \Sigma(\pi_{\tilde\xi, \tilde\eta}) +
\frac{\tau_1 + \theta -t}{t} \int \log x \ d\pi_{\tilde\xi,
\tilde\eta}(x) +  \frac{\tau_2 -t}{t} \int \log(1-x) \
d\pi_{\tilde\xi, \tilde\eta}(x)  \\=  B\Big(\frac{\tau_1
+\theta-t}{t}, \frac{\tau_2 -t}{t}\Big) \ee and then, on ${\mathcal
A}(\xi T^{-1})$ the infimum is uniquely realized in $\pi_{\tilde\xi,
\tilde\eta}$ and its value is
\[
\theta\xi +t^2 \left[B\Big(\frac{\tau_1 -t}{t}, \frac{\tau_2
-t}{t}\Big) - B\Big(\frac{\tau_1 +\theta-t}{t}, \frac{\tau_2
-t}{t}\Big)\right]\,.
\]
Finally a small computation leads to (\ref{idetJ}) and
(\ref{defxiJ}).
\medskip

\noindent{\it{Proof of Lemma \ref{min}}}
 In \cite{Saff} p.241, it is proved that the infimum of
$$ \int\!\int -\log |x-y| d\mu(x) d\mu(y) -2 \zeta_1\int \log (1-x)\!\ d\mu(x)-2\zeta_2\int \log (1+x)\!\  d\mu(x)  $$
among the probability measures $\mu$ on $[-1, +1]$ is achieved by
$$d\mu (y) = K(b_-,b_+) \frac{\sqrt{(y - b_-)(b_+ - y)}}{2\pi (1-y^2)}\!\  \mathbf{1}_{[b_-,b_+]}(y)\ dy\,,$$
where $b_\pm = \theta_2^2 -\theta_1^2 \pm \sqrt\Delta$ with \ben
\nonumber \label{theta} \theta_i =
\frac{\zeta_i}{1+\zeta_1+\zeta_2}\ , \ i=1,2\ &,& \ \Delta =
\left[1-(\theta_1+\theta_2)^2
\right]\left[1-(\theta_1-\theta_2)^2\right]\,. \een and $K(b_-,b_+)$
is a normalizing constant. With the push-forward by the function $x
\rightarrow (x+1)/2$, we get the result. $\Box$

\section{Appendix 1 : Some properties of $\ell = \log \Gamma$ and $\Psi$}
\label{appendix} From the Binet formula (
\cite{astig} pp. 258-259 or  \cite{bateman} p.21), we have \ben
\label{bin1}
\ell (x) &=& (x-\frac{1}{2})\log x -x +1 + \int_0^\infty f(s)[e^{-sx} - e^{-s}]\ \! ds\\
\label{bin2}
&=& (x-\frac{1}{2})\log x -x + \frac{1}{2} \log(2\pi) + \int_0^\infty f(s)e^{-sx}\ \! ds\,.
\een
where the function $f$ is defined by
\ben
\label{propf}
 f(s) = \left[\frac{1}{2}-\frac{1}{s}+ \frac{1}{e^s -1}\right]\frac{1}{s} =
 2\sum_{k=1}^\infty \frac{1}{s^2 + 4\pi^2 k^2}\,,
\een and satisfies for every $s \geq 0$ \ben \label{proprf} 0 < f(s)
\leq
\label{bina} f(0)= 1/12 \ , \  \ 0 < \left(sf(s)+ \frac{1}{2}\right) < 1\,.
\een By differentiation
\ben
\label{bin3}
\log x - \Psi (x) = \frac{1}{2x} +\int_0 ^\infty s f(s) e^{-sx}\, ds = \int_0 ^\infty e^{-sx} \left(sf(s)+ \frac{1}{2}\right)
 ds\,.
\een
As easy consequences, we have, for every $x > 0$
\ben
\label{supx}
0 &<& x\left(\log x - \Psi (x)\right) \leq 1\,,\\
\label{supx2} 0 &<& \log x - \Psi(x) - \frac{1}{2x} \leq
\frac{1}{12x^2}\,. \een Differentiating again we see that for $q\geq
1$ \ben \label{polygamma} \Psi^{(q)}(z) = (-1)^{q-1} q! z ^{-q} +
(-1)^{q-1} \int_0^\infty e^{-sz} s^q  \left(sf(s)+
\frac{1}{2}\right)\ ds \een and then \ben \label{restepsi}
 |\Psi^{(q)}(z) - (-1)^{q-1} q! z ^{-q}| \leq z^{-q-1} q!\,.
 \een
\section{Appendix 2 : Identification
 of the McKay distribution}
\label{appendixMK}
The reader is recalled that, for $u'$ and $v'$ positive numbers\footnote{we use the
symbol   $v'$ (hence $u'$) not to confuse with $\beta$ already
defined.} such that $u' + v' > 1$, 
\cite{capitaine} defined the probability measure \be 
\hbox{CC}_{u', v'} := (1- u')^+\delta_0 + (1-v')^+ \delta_1 +\left[1
- (1- u')^+ - (1-v')^+ \right]\pi_{a_-,a_+}\,, \ee where \be
 (a_-,a_+) = a_\pm\Big(\frac{u'}{u' + v'} , 1 -
\frac{1}{u' + v'}\Big)\,. \ee We present now three identifications
of this distribution connected  with free probability.

 For $k \not=
0$, let $D_k$ the dilatation operator by factor $k$. For $p \leq 1$,
let $\b_p$ denote the Bernoulli distribution of parameter $p$. At
last, let $\boxplus$ (resp. $\boxtimes$) denote the additive (resp.
multiplicative) free convolution.
\medskip

\noindent 1) Rewriting the distribution with the notation of 
\cite{Dem}, we find four cases

$\bullet$ Situation $I$ : $\min(u' , v') \geq 1$, no Dirac mass,
$$\sigma_- = \frac{u'}{u'+v'} \ , \ \sigma_+ = 1 -
\frac{1}{u'+v'} \ , \ u' = \frac{\sigma_-}{1-\sigma_+} \ , \ v' =
\frac{1-\sigma_-}{1-\sigma_+}$$
\be
CC_{u' , v'} = \pi_{a_-, a_+} = D_{1-\sigma_+}
(\b_{\sigma_-})^{\boxplus\frac{1}{1- \sigma_+}} \ee

$\bullet$ Situation $II$ : $u' < 1 \leq  v'$, one Dirac mass at $0$
$$\sigma_- = \frac{1}{u'+v'} \ , \ \sigma_+ = 1 - \frac{u'}{u'+v'} \ , \ u' = \frac{1 - \sigma_+}{\sigma_-} \ , \ v' =  \frac{\sigma_+}{\sigma_-}$$
\be
CC_{u' , v'} &=& (1- u')\delta_0 + u'\pi_{a_-,a_+}\\
 &=& D_{\sigma_-} (\b_{1 - \sigma_+})^{\boxplus\frac{1}{\sigma_-}}
\ee

$\bullet$ Situation $III$ : $v' < 1 \leq  u'$, one Dirac mass at $1$
$$\sigma_- = 1 - \frac{1}{u'+v'} \ , \ \sigma_+ = \frac{u'}{u'+v'} \ , \ u' = \frac{\sigma_+}{1 - \sigma_-} \ , \ v' =  \frac{1- \sigma_+}{1- \sigma_-}$$
\be
CC_{u' , v'} &=& (1-v') \delta_1  + v' \pi_{a_-,a_+}\\
 &=& D_{1- \sigma_-} (\b_{\sigma_+})^{\boxplus\frac{1}{1-\sigma_-}}
\ee

$\bullet$ Situation $IV$ : $\max(u', v') < 1$, two Dirac masses (at
$0$ and at $1$)
$$\sigma_- = 1 - \frac{u'}{u'+v'} \ , \ \sigma_+ = \frac{1}{u'+v'} \ , \ u' = \frac{1- \sigma_-}{\sigma_+} \ , \ v' =  \frac{\sigma_-}{\sigma_+}$$
\be
CC_{u' , v'} &=& (1- u')\delta_0 + (1-v') \delta_1 + (u' + v' -1)\pi_{a_-,a_+}\\
&=& D_{\sigma_+} (\b_{1-\sigma_-})^{\boxplus\frac{1}{\sigma_+}}\,.
\ee
\medskip

\noindent 2) There is a  connection with the  family of free Meixner
law (\cite{bryc1}, \cite{bryc2}, \cite{bryc3}). Indeed,
computing the mean $m$ and variance $V$ of the distribution
$CC_{u',v'}$, we get
\medskip

\begin{center}
\begin{tabular}{|c|c|c|}
\hline
Situation & m & V \\
\hline
I & $\sigma_-$ &  $\sigma_- (1-\sigma_-) (1-\sigma_+)$\\
\hline
II & $1-\sigma_+$  & $\sigma_- \sigma_+(1-\sigma_+)$ \\
\hline
III &  $\sigma_+$ & $(1-\sigma_-) \sigma_+ (1-\sigma_+)$\\
\hline
IV & $1- \sigma_-$ & $\sigma_- \sigma_+(1-\sigma_-)$\\
\hline
\end{tabular}
\end{center}
\medskip

so that, in all cases \be m  = \frac{u'}{u' + v'}\ , \ V  =
\frac{u'v'}{(u' + v')^3}\,. \ee We see that fixing $u' + v' =
s^{-1}$, we get $V = s^2 m (1-m)$, and then up to an affine
transformation we find the "free binomial type law" as in \cite{bryc3}
example vi p.18 or \cite{bryc1} example 6 p.8. It could also be seen
starting from the above formulae using dilatations and free
convolutions and comparing with formula (7) page 6 in \cite{bryc1}.
\medskip

\noindent 3) Finally, we quote the correspondence with the results
of  \cite{Collins} who claimed that for $0 < p_- < p_+ < 1$
\be
\b_{p_-}\boxtimes\b_{p_+} = (1 - p_-) \delta_0 +    (p_- + p_+ -
1)^+ \delta_1 + C_{a_- , a_+}^{-1} \pi_{a_- , a_+} \ee where $a_\pm
= a_\pm(1-p_- , p_+)$. In \cite{hiai2}, formula (2.8) the authors
consider the same distribution.

$\bullet$ Situation $II$ : $p_- + p_+ - 1 < 0$,  $\sigma_- = p_+$ , $\sigma_+ = 1 - p_-$
\be\b_{p_-}\boxtimes\b_{p_+} &=& \sigma_+ \delta_0 + (1 - \sigma_+) \pi_{a_- , a_+}\\
&=& \sigma_- CC_{u' , v'} + (1 - \sigma_-) \delta_0 \ee

$\bullet$ Situation $IV$ : $p_- + p_+ - 1 > 0$, $\sigma_- = 1 - p_-$ , $\sigma_+ = p_+$
\be
\b_{p_-}\boxtimes\b_{p_+}  &=& \sigma_- \delta_0 + (\sigma_+ - \sigma_-)\delta_1 + (1 - \sigma_+)\pi_{a_- , a_+}\\
&=& \sigma_+ CC_{u' , v'} + (1 - \sigma_+) \delta_0\,. \ \ \ \ \Box \ee
\bigskip

{\bf Acknowledgement}
\medskip

I thank Catherine Donati-Martin for introducing me to some references. I am also grateful to the anonymous referee for a careful reading of the manuscript.

\end{document}